\documentclass[11pt]{article}

\usepackage[T1]{fontenc}
\usepackage[utf8]{inputenc}

\usepackage{amsmath}
\usepackage{amsthm}
\usepackage{mathtools}
\usepackage{amsfonts}
\usepackage{url}
\usepackage{microtype}

\usepackage{lmodern}

\usepackage[top=1.1in, bottom=1.1in, left=1.25in, right=1.25in]{geometry}

\usepackage[font=small]{caption}
\usepackage[font=footnotesize]{subcaption}

\usepackage{graphicx}
\usepackage{afterpage}
\usepackage{setspace}
\usepackage{xcolor}
\usepackage{color}

\usepackage[sf,bf,medium]{titlesec}

\usepackage[numbers,sort&compress]{natbib}

\DeclareMathOperator{\Res}{Res}

\DeclareMathOperator{\realpart}{Re}

\usepackage{isomath}
\renewcommand{\phi}{\varphi}

\usepackage{bm}
\renewcommand{\vec}{\vectorsym}
\newcommand{\vct}{\vectorsym}

\newcommand{\surfdiv}{\nabla_\Gamma \cdot}
\newcommand{\surfdivp}{\nabla_\Gamma' \cdot}
\newcommand{\surfdivl}{\nabla_{\Gamma_\ell} \cdot}
\newcommand{\surfgrad}{\nabla_\Gamma}
\newcommand{\surfgradl}{\nabla_{\Gamma_\ell}}
\newcommand{\surflap}{\Delta_\Gamma}
\newcommand{\surflapm}{\Delta_{\Gamma,m}}
\newcommand{\surflapl}{\Delta_{\Gamma_\ell}}

\newcommand{\bbR}{\mathbb R}

\newcommand{\bn}{\vct{ n}}

\newcommand{\bx}{\vct{ x}}

\newcommand{\bL}{\vct{ L}}

\newcommand{\bE}{\vectorsym{E}}
\newcommand{\bcB}{\bm{{\mathcal B}}}
\newcommand{\bcD}{\bm{{\mathcal D}}}
\newcommand{\bcE}{\bm{{\mathcal E}}}
\newcommand{\bcA}{\bm{\mathcal A}}
\newcommand{\bcH}{\bm{\mathcal H}}
\newcommand{\bcJ}{\bm{\mathcal J}}
\newcommand{\bff}{\vct{ f}}
\newcommand{\bF}{\vct{ F}}
\newcommand{\bH}{\vct{ H}}
\newcommand{\bJ}{\vct{ J}}
\newcommand{\bK}{\vct{ K}}

\newcommand{\cI}{\mathcal I}

\newcommand{\cS}{\mathcal S}

\newcommand{\cC}{\mathcal C}
\newcommand{\cU}{\mathcal U}

\newtheorem{remark}{\sffamily Remark}

\newcommand{\cK}{\mathcal K}

\newcommand\beqal{\begin{equation}\begin{aligned}}
\newcommand\eeqal{\end{equation}\end{aligned}}

\newcommand\bzero{\vct{ 0}}

\newcommand\btheta{\vct{ \theta}}

\newcommand\bsf{\vct{ f}}
\newcommand\ihat{\vct{{i}}}
\newcommand\jhat{\vct{{j}}}
\newcommand\khat{\vct{{k}}}
\newcommand\rhat{\vct{{r}}}
\newcommand\tauhat{\vct{{ \tau}}}
\newcommand\thetahat{{\vct{{\theta}}}}
\newcommand\nhat{\vct{ {n}}}

\newcommand\bl{\vct{ l}}
\newcommand\bu{\vct{ u}}

\newcommand\bQ{\vct{ Q}}

\newcommand\cO{\mathcal O}

\newcommand\btau{\vct{\tau}}
\newcommand\bA{\vct{A}}
\newcommand\ba{\vct{a}}
\newcommand\bh{\vct{h}}

\newcommand\bEtot{\vct{E}^\textnormal{tot}}
\newcommand\bHtot{\vct{H}^\textnormal{tot}}
\newcommand\bEin{\vct{E}^\textnormal{in}}
\newcommand\bHin{\vct{H}^\textnormal{in}}

\newcommand\Ein{E^\textnormal{in}}
\newcommand\Hin{H^\textnormal{in}}

\title{\bf\sffamily A high-order wideband direct 
  solver for electromagnetic
  scattering from bodies of revolution}

\numberwithin{equation}{section}

\author{Charles L. Epstein\footnote{Mathematics Department, University
    of Pennsylvania, Philadelphia, PA. Email: {\tt
      cle@math.upenn.edu}.  Research supported in part by NSF award
     DMS-1507396.},
  \ Leslie Greengard\footnote{Courant Institute, New York University,
    New York, NY and Flatiron Institute, New York, NY.  Email: {\tt
      greengard@cims.nyu.edu}.  Research supported in part by
    the Office of the Assistant Secretary of Defense for Research and
    Engineering and AFOSR under NSSEFF Program award
    FA9550-10-1-0180.}, \ and Michael O'Neil~\footnote{Courant
    Institute, New York University, New York, NY. Email: {\tt
      oneil@cims.nyu.edu}.  Research supported in part by the
    Office of the Assistant Secretary of Defense for Research and
    Engineering and AFOSR under NSSEFF Program award FA9550-10-1-0180
    and the Office of Naval Research under awards N00014-17-1-2059 and
    N00014-17-1-2451.}}

\date{\today}

\usepackage{amsopn}

\begin{document}

\maketitle

\begin{abstract}
  The generalized Debye source representation of time-harmonic
  electromagnetic fields yields well-conditioned second-kind integral
  equations for a variety of boundary value problems, including
  the problems of scattering from perfect
  electric conductors and dielectric bodies.  Furthermore, these
  representations, and resulting integral equations, are fully stable
  in the static limit as $\omega \to 0$ in multiply connected
  geometries. In this paper, we present the
  first high-order accurate solver based on this representation for
  bodies of revolution.
  The resulting solver uses a Nystr\"om discretization of a one-dimensional
  generating curve and high-order integral equation methods for
   applying and inverting surface differentials.  The accuracy and
  speed of the solvers are demonstrated in several numerical examples.

  \vspace{.5\baselineskip}
\end{abstract}




\section{Introduction}
\label{sec_intro}

In isotropic homogenous media, the time-harmonic
Maxwell's equations (THME) governing the
propagation of electric and magnetic fields are given by~\cite{papas}
\begin{equation}\label{eq_maxwell}
  \begin{aligned}
    \nabla \times \bcE &= i \omega \mu \bcH, 
    &\qquad \nabla \times \bcH &= \bcJ - i\epsilon\omega \bcE,\\
    \nabla \cdot \bcE &= \frac{\rho}{\epsilon}, & \nabla \cdot \bcH &= 0,
  \end{aligned}
\end{equation}
where $\bcE$ and $\bcH$ are the electric and magnetic fields,
respectively, and $\omega$ is the angular frequency. The electric
charge is denoted by $\rho$ and the electric current by $\bcJ$.
The electric
permittivity and magnetic permeability are given by $\epsilon$ and
$\mu$, respectively, and assumed to be piecewise constant in all that
follows.
 It will be useful to define the wavenumber $k$ as
$k = \omega\sqrt{\epsilon\mu}$, and we assume that $\Im (k) \geq 0$ as well.
The implicit
time dependence of $e^{-i\omega t}$ on $\bcE$, $\bcH$ has been suppressed.

Maxwell's equations in the above form provide a surprisingly accurate
approximation to physical phenomena, and therefore their accurate
numerical solution is of the utmost importance in many fields, such as
radio-frequency component modeling~\cite{kapur-1998}, meta-material
simulation~\cite{gimbutas-2013}, and radar cross-section
characterizations~\cite{youssef1989}.  Several computational methods
for the solution to~\eqref{eq_maxwell} are used in practice, including
finite difference, finite element, spectral methods, and integral
equation
methods~\cite{monk1992finite,kirsch1995finite,nedelec,yee_1966}.
Finite difference and finite element methods are particularly useful
when the material parameters vary in space, especially when that
variation is anisotropic.

On the other hand, when applicable -- namely in complex geometries
with piecewise constant material parameters -- integral
equation methods have significant advantages. They reduce the
dimensionality of the problem to the boundary alone, satisfy the 
necessary radiation conditions in exterior domains, do not suffer from
grid-based dispersion errors (particularly problematic at high
frequencies),
and are compatible with a variety of fast
algorithms, leading to asymptotically optimal schemes.
Two standard problems in electromagnetics -- scattering from perfect
electric conductors and dielectric bodies -- fit exactly these
assumptions~\cite{jackson}.  
While high-order accurate
methods in arbitrary geometries in three dimensions are still under
active development,
there are several important applications of Maxwell's equations in
scattering from bodies of revolution~\cite{Chewfastbook,mautz1979electromagnetic,
gedney1990use,Kucharski2000,helsing2015determination,liu2016efficient,viola},
where surface representation and quadrature are more straightforward.
In this context, it is possible to build very high-order accurate, quasi-optimal
solvers using separation of variables and the Fast Fourier Transform
(FFT) which scale roughly as~$\cO(N^{3/2})$, where $N$ is the total number
of unknowns needed to fully discretize the body's surface.

Despite the applicability of integral
equation-based solvers to bodies of revolution, however, wideband solvers have
yet to be implemented that are robust in the genus-one case and 
in the low-frequency regime. There are two principal difficulties here:
(1) constructing high-order accurate Nystr\"om methods
(including quadrature) for the so-called modal Green's
functions~\cite{conway_cohl} is somewhat technical, and (2) existing
integral equation methods suffered from ill-conditioning due to
\emph{dense-mesh breakdown} and \emph{topological
  breakdown}~\cite{cools-2009,EpGr,EpGrOn} in the static limit.  The
algorithm of the present paper, based on the recent generalized Debye source
formulations of Maxwell fields~\cite{EpGr,EpGrOn}, overcomes both of
these obstacles.  

The past few years have seen several methodological developments in
electrostatic and acoustic scattering from bodies of revolution that
should be noted~\cite{helsing_2014,young}.  The vector-valued case is
treated in~\cite{liu2016efficient}, but with the method of fundamental
solutions~\cite{barnett2008stability,fairweather1998method} rather
than an on-surface boundary integral discretization.  In this work, we
provide a thorough discussion of the calculations required in the
vector case, as it is more involved than the scalar case.  There is a
substantial literature on scattering from bodies of revolution in the
electrical engineering community which we do not seek to review here.
We refer the reader to the papers and books mentioned above. We will
point out relevant previous work when discussing particular topics:
quadrature, kernel evaluation, and integral equation conditioning.

The paper is organized as follows: Section~\ref{sec_gendebye}
introduces the generalized Debye source formulation of Maxwell fields
on general surfaces and formulates the relevant boundary value
problems for scattering from perfect conductors and dielectric
(penetrable) bodies.  Section~\ref{sec_revolution} derives the
necessary formalism required for separating vector-valued functions in
cylindrical coordinates on surfaces of revolution.
Section~\ref{sec_solver} describes a high-order Nystr\"om-type solver
for electromagnetic scattering. This includes a discussion of
quadrature, kernel evaluation, and surface differential application.
We provide several numerical examples demonstrating our algorithm in
Section~\ref{sec_examples} and then discuss extensions and drawbacks
of the scheme in Section~\ref{sec_conclusions}.

\section{Generalized Debye sources}
\label{sec_gendebye}

Classically, most integral equations for electromagnetic scattering
problems have been derived in the Lorenz gauge. More precisely,
the electric and magnetic fields are represented
using a scalar potential $\varphi$ and vector potential
$\bcA$~\cite{papas}:
\begin{equation}
  \begin{aligned}
    \bcE &= i\omega \bcA - \nabla \varphi \\
    \bcH &=  \nabla \times \bcA,
  \end{aligned}
\end{equation}
with
\begin{equation}\label{eq_pots}
  \bcA(\bx) = \frac{1}{\mu} \int \frac{e^{ik|\bx-\bx'|}}{4\pi | \bx-\bx'
    |} \, \bcJ(\bx') \, d\bx', \qquad
  \varphi(\bx) = \frac{1}{\epsilon} \int \frac{e^{ik|\bx-\bx'|}}{4\pi | \bx-\bx'
    |} \, \frac{\nabla \cdot \bcJ(\bx') }{i\omega} \, d\bx',
\end{equation}
so that $\nabla \cdot \bcA = i\omega\epsilon\mu\varphi$.  In the case
of the perfect electric conductor, using these representations yields
the Electric Field Integral Equation (EFIE) by enforcing the condition
that the tangential electric field vanishes. The Magnetic Field
Integral Equation (MFIE) is derived from the fact that the tangential
magnetic field must equal the surface current $\bcJ$. By taking linear
combinations of these equations, one can derive various Combined Field
Integral Equations (CFIE)~\cite{colton_kress} to avoid spurious
resonances.  However, unless additional preconditioners are
used~\cite{contopanagos-2002}, the EFIE and CFIE suffer from what is
known as \emph{dense-mesh low-frequency breakdown} due to the division
by $\omega$ in the scalar potential term, leading to catastrophic
cancellation.  This can be circumvented by introducing electric charge
as an additional variable, subject to the consistency condition
$\nabla \cdot \bcJ = i\omega \rho$, which is implied by Maxwell's
equations~\eqref{eq_maxwell}. In that case, we define the scalar
potential as $\varphi = \cS_k \rho/\epsilon$. Here, $\cS_k$ denotes
the single-layer potential operator, as in~\eqref{eq_pots}.
Additional boundary conditions must be added because of the extra
unknown, but this is
straightforward. (See~\cite{taskinen-2006,vico_2016,vico-2013} for
various \emph{charge-current formulations} of scattering problems.)
However, as noted in the introduction, the second form of
low-frequency breakdown occurs in Maxwell's equations, which we refer
to as \emph{topological low-frequency breakdown}~\cite{EpGr}.  This
form of breakdown is characterized by increased ill-conditioning as
$\omega \to 0$, with a singular limit in the static case.  The
generalized Debye source representation of Maxwellian fields overcomes
both, by taking the genus into account and by using scalar variables
rather than currents to avoid catastrophic cancellation~\cite{EpGr,
  EpGrOn}.  (An alternative robust formulation is the Decoupled
Potential Integral Equation (DPIE)~\cite{vico_2016}, which results in
separate integral equations for $\bcA$ and $\varphi$.)

For charge- and current-free regions with piecewise constant material
parameters $\epsilon$, $\mu$, a simple re-scaling of the
equations~\eqref{eq_maxwell} is useful to reduce the dependency to
only a single parameter~$k$. To this end, we define
\begin{equation}
  \bEtot = \sqrt{\epsilon} \bcE, \qquad
  \bHtot = \sqrt{\mu} \bcH,
\end{equation}
The THME then take the form:
\begin{equation}
\label{eq_maxwell_k}
  \begin{aligned}
    \nabla \times \bEtot &= i k \bHtot, 
    &\qquad \nabla \times \bHtot &= - ik \bEtot,\\
    \nabla \cdot \bEtot &= 0, & \nabla \cdot \bHtot &= 0.
  \end{aligned}
\end{equation}
where, as before, $k=\omega \sqrt{\epsilon\mu}$.  In electromagnetic
scattering problems, the \emph{total} fields $\bEtot$,~$\bHtot$
consist of the superposition of a known, incoming field~$\bEin$,
$\bHin$ and an unknown scattered field $\bE$, $\bH$:
\begin{equation}
  \begin{aligned}
    \bEtot &= \bE + \bEin \\
    \bHtot &= \bH + \bHin.
  \end{aligned}
\end{equation}
Boundary conditions are enforced on the total fields, and the known
incoming fields provide data for the resulting integral equations.

\begin{figure}[t]
  \centering
  \includegraphics[width=.6\linewidth]{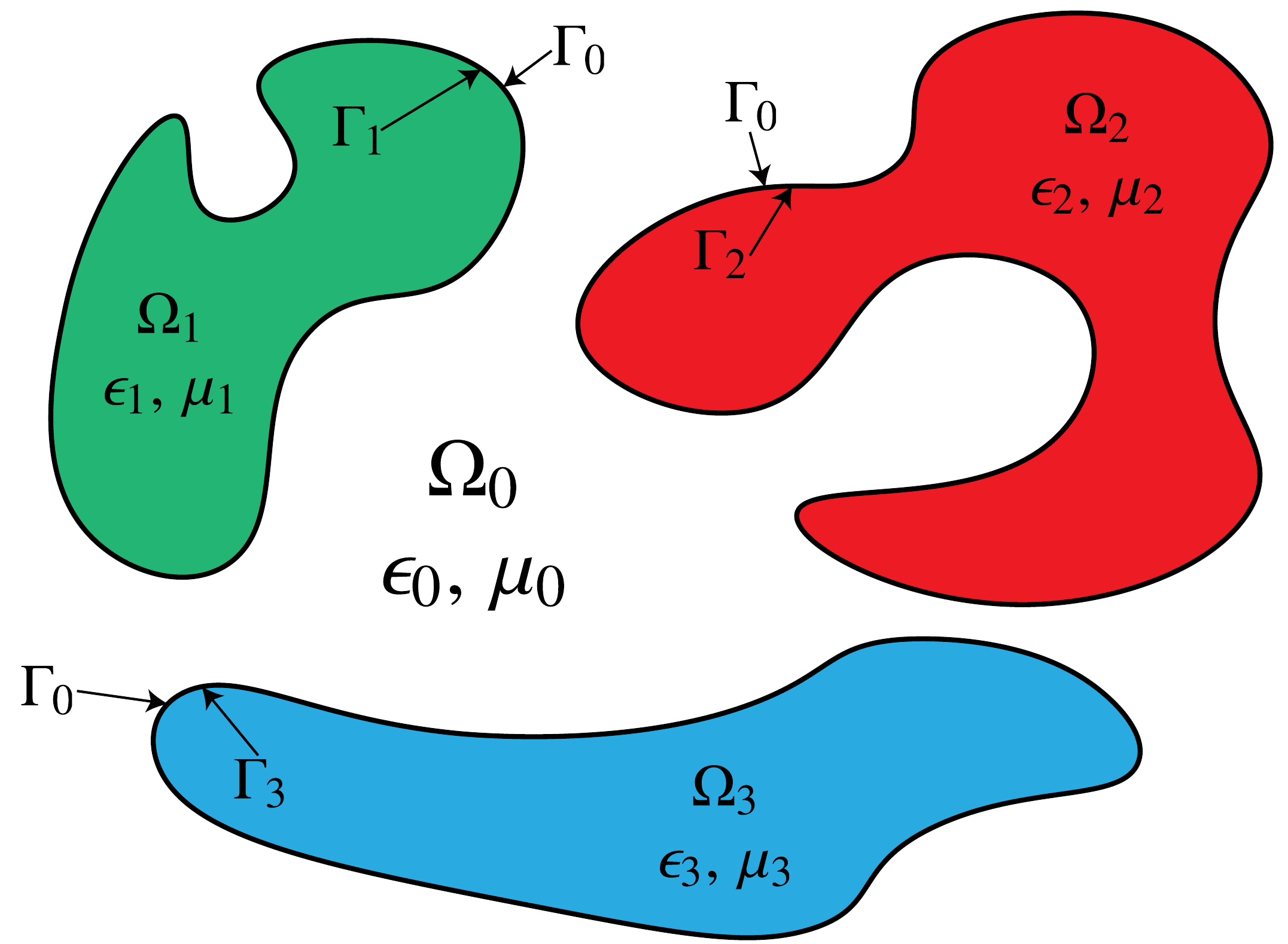}
  \caption{Boundary, domain, and material specifications for
    piecewise-constant geometries.}
  \label{fig_dielectric}
\end{figure}

In the case of
multiple-material (non-conducting) penetrable scattering, the values
of $\epsilon$, $\mu$ in region $\ell$ will be denoted by
$\epsilon_\ell$, $\mu_\ell$, and all relevant quantities will inherit
the subscripts as well (the electromagnetic fields, scalar/vector
functions supported on interfaces, etc.). See
Figure~\ref{fig_dielectric}. The parameters for the
unbounded region (air or vacuum) will be denoted by $\epsilon_0$,
$\mu_0$.  In the case of exterior scattering from a collection of perfectly
conducting objects, there is only one region (the exterior)
with non-zero fields, and the subscript $\ell$ can be suppressed.
We now briefly describe a representation of electromagnetic fields using
generalized Debye sources which is compatible with the multiple region
configuration. 

The generalized Debye source representation makes use of \emph{non-physical} 
variables - using both potentials ($\bA, \varphi$) and \emph{anti-potentials}
($\bQ, \psi$):
\begin{equation}\label{eq_gendebye}
  \begin{aligned}
    \bE_\ell &= ik_\ell \bA_\ell - \nabla \varphi_\ell - \nabla \times \bQ_\ell, \\
    \bH_\ell &= ik_\ell \bQ_\ell - \nabla \psi_\ell + \nabla \times \bA_\ell,
  \end{aligned}
\end{equation}
where the potentials are given by
\begin{equation}
  \begin{aligned}
    \bA_\ell &=  \mathcal S_{k_\ell} \bJ_\ell, \qquad& \varphi_\ell &=
    \cS_{k_\ell} \rho_\ell, \\
    \bQ_\ell &=  \mathcal S_{k_\ell} \bK_\ell , & \psi_\ell &=
    \cS_{k_\ell} \sigma_\ell,
  \end{aligned}
\end{equation}
and the layer potentials $S_{k_\ell}$ are defined on the support of their
argument. Thus,
\begin{equation}
  \varphi_\ell(\bx) = \int_{\Gamma_\ell} \frac{e^{ik|\bx-\bx'|}}
  {4\pi |\bx-\bx'|} \rho_\ell(\bx') \, da(\bx'),
\end{equation}
with analogous expressions defining the other potentials in terms of their
vector or scalar source densities ($\bJ_\ell, \bK_\ell$, and $\sigma_\ell$).
Here, in a slight abuse of notation, $\rho$ is 
not the true electric charge, but a charge-like 
non-physical quantity.
The tangential vector fields $\bJ_\ell$, $\bK_\ell$ and the scalar
functions $\rho_\ell$, $\sigma_\ell$ are defined on the boundary of
the $\ell$th region~$\Gamma_\ell$ (see Figure~\ref{fig_dielectric}),
and are responsible for representing the fields
in the $\ell$th region.
In order for Maxwell's equations to be satisfied in a
source-free region, the potentials must satisfy the relations
\mbox{$\nabla \cdot \bA_\ell = ik_\ell \phi_\ell$} and \mbox{$\nabla \cdot
\bQ_\ell = ik_\ell \psi_\ell$}. It is straightforward to show that this 
holds when the currents and charges satisfy the 
consistency conditions:
\begin{equation}
  \surfdivl \bJ_\ell
  = ik_\ell \rho_\ell, \qquad \surfdivl \bK_\ell = ik_\ell \sigma_\ell.
\label{consistcond}
\end{equation}

Rather than using these conditions to \emph{define} 
$\rho_\ell, \sigma_\ell$, leading to catastrophic cancellation,
or adding them as additional unknowns, requiring their imposition as
additional constraints, we will construct $\bJ_\ell$ and $\bK_\ell$
in such a way that \eqref{consistcond} is automatically satisfied.
That construction lies at the heart of the generalized Debye formalism,
and the scalar surface densities $\rho_\ell$, $\sigma_\ell$ are referred
to as \emph{generalized Debye sources}. 
%
Thus, we define the currents $\bJ_\ell$ and $\bK_\ell$ in terms of 
an explicit Hodge decomposition into
divergence-free, curl-free, and harmonic components
along each $\Gamma_\ell$:
\begin{equation}\label{eq_hodge}
  \begin{aligned}
    \bJ_\ell &= \surfgradl \alpha_{J_\ell} + \nhat \times \surfgradl
        \beta_{J_\ell} + \bJ^h_\ell, \\
    \bK_\ell &= \surfgradl \alpha_{K_\ell} + \nhat \times \surfgradl \beta_{K_\ell} + \bK^h_\ell,
  \end{aligned}
\end{equation}
where 
\begin{equation}\label{surflapeq}
\surflapl \alpha_{J_\ell} = ik_\ell \rho_\ell, \qquad 
\surflapl \alpha_{K_\ell} = ik_\ell \sigma_\ell,
\end{equation}
$\surfgradl$ denotes the surface gradient,
$\surflapl$ denotes the surface Laplacian, 
and $\bJ^h_\ell$, $\bK^h_\ell$ are harmonic vector fields on $\Gamma_\ell$.
Recall that a tangential vector field $\bL$ is {\em harmonic} if it satisfies 
the conditions:
\begin{equation}
\surfdivl \bL = 0, \qquad \surfdivl \nhat \times \bL = 0.
\end{equation}
The dimension of the linear space of harmonic vector fields on each
boundary component~$\Gamma_\ell$ is simply $2g_\ell$,
where $g_\ell$ is the genus of (the smooth) boundary $\Gamma_\ell$.
Thus, we will write
\begin{equation}
  \bJ^h_\ell = \sum_{j = 1}^{2g_\ell} a_{\ell j} \, \bh_{\ell j}, \qquad 
  \bK^h_\ell = \sum_{j = 1}^{2g_\ell} b_{\ell j} \, \bh_{\ell j},
\end{equation}
where the $\bh_{\ell j}$ form a basis for the harmonic vector fields.
It is the coefficients $a_{\ell j}$, $b_{\ell j}$ that must be determined
based on boundary data.
It is important to note that the equations~\eqref{surflapeq}
guarantee that~\eqref{consistcond} is satisfied.
One slight complication, addressed below, is that the
sources $\rho_\ell$, $\sigma_\ell$ are required to be
mean-zero, since the surface Laplacian is invertible only in the space
of mean-zero functions~\cite{EpGr}.
Finally, $\beta_{J_\ell}$ and $\beta_{K_\ell}$ are chosen to be 
linear in $\rho_\ell$ and $\sigma_\ell$ in a manner that depends on the precise
boundary conditions of interest. 

\begin{figure}[t]
  \centering
  \includegraphics[width=.4\linewidth]{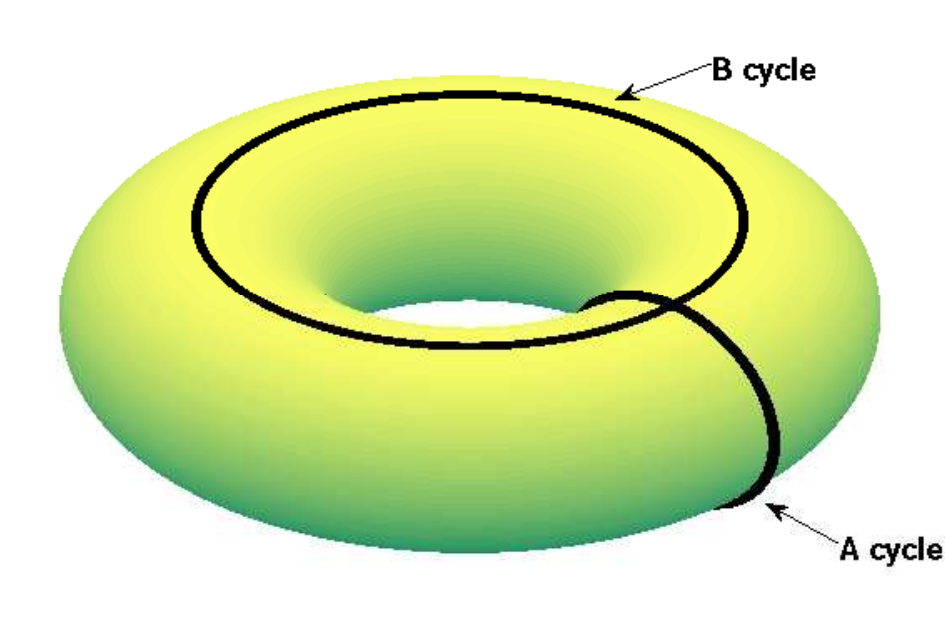}
  \caption{Loops spanning the interior and exterior homologies of a
      genus 1 surface.}
  \label{fig_cycles}
\end{figure}

For boundary-value problems in multiply-connected geometries, the generalized
Debye approach requires a set of conditions to fix the harmonic components
beyond the local boundary conditions imposed on the field quantities
$\bE$, $\bH$.
These (non-local) conditions are determined by considering $n_g = 2\sum_\ell
g_\ell$ separate
line integrals of $\bEtot$, $\bHtot$ along loops contained in $\Gamma_\ell$ 
which enclose \emph{interior and exterior holes} -- that is to say,
these loops \mbox{$C_1$, \ldots, $C_{n_g}$} must span the first
homology group of $\Gamma_\ell.$ 
They take the general form:
\begin{equation}
\int_{C_j} \bE \cdot d\bl = f_j,\quad
\int_{C_j} \bH \cdot d\bl = g_j,
\end{equation}
for $j = 1, \ldots, n_g$ and some set of constants $f_j,g_j$. In scattering
problems, the
constants~$f_j, g_j$ usually correspond to some functional of the
incoming fields and depend, again, on the particular boundary conditions.
We will specify the precise constraints below for perfect conductors and 
dielectrics.
See Figure~\ref{fig_cycles} for an illustration of a possible set of 
loops $C_j$ in the genus 1 case.

\begin{remark}
Since the generalized Debye construction enforces that 
the electric and magnetic fields satisfy the Maxwell equations, it
can be used in applications which do not arise as
\emph{scattering problems}~\cite{EpGrOn2}. This is not true of
many combined-field and charge-current representations, which work
in the scattering setting only because the right-hand side is
not arbitrary - it is derived from an incoming Maxwellian electromagnetic field
(see, for example, \cite{vico-2013}).
\end{remark}

We now turn to the particular representation used for electromagnetic
fields in the exterior of perfect conducting objects.

\subsection{The perfect electric conductor (PEC)}
\label{sec_pec}

We describe the generalized Debye representation for a single PEC
object, since the extension to multiple 
objects is straightforward.  Along the surface~$\Gamma$
of $\Omega$ (a PEC), it is well known that 
\begin{equation}\label{eq_pecbcs}
    \nhat \times \bcE = \bzero, \qquad \nhat \cdot \bcH = 0.
\end{equation}
In the case of scattering from a PEC, the current-like vector
fields~$\bJ$ and~$\bK$
in the generalized Debye source representation can be constructed
as~\cite{EpGr}:
\begin{equation}
  \label{eq_jbuild}
  \begin{aligned}
    \bJ &= ik \left( \surfgrad \surflap^{-1} \rho - \nhat \times 
      \surfgrad \surflap^{-1} \sigma \right) + \bJ^h, \\
    \bK &= \nhat \times \bJ.
  \end{aligned}
\end{equation}
Here, 
$\surflap^{-1}$ denotes the inverse of the Laplace-Beltrami
operator (the surface Laplacian) along $\Gamma$ restricted to mean-zero functions. 
On the space of mean-zero functions defined along~$\Gamma$, 
this operator is uniquely invertible. Using
this construction for $\bJ$ and $\bK$ and enforcing the
preconditioned~\cite{nedelec} homogeneous boundary conditions
\begin{equation}\label{eq_pecinteq1}
  \begin{aligned}
    \cS_0 \surfdiv \bEtot &= 0, &\quad &\text{on } \Gamma, \\
     \nhat \cdot \bHtot &= 0,  & &\text{on } \Gamma, \\
     \int_{C_j} \bEtot \cdot d\bl &= 0, & &\text{for } j = 1,\ldots,
     2g
  \end{aligned}
\end{equation}
yields a second-kind integral equation for $\rho$, $\sigma$, and
$\bJ^h = \sum_{j = 1}^{2g_\ell} a_{j} \, \bh_{j}$
which is uniquely invertible for all values of $k$ with
non-negative imaginary part~\cite{EpGr}. By construction, the
scattered fields automatically satisfy the Silver-M\"uller radiation
conditions at infinity.
Here $\cS_0$ is the
zero-frequency single-layer potential operator.
The system in~\eqref{eq_pecinteq1} can be written more explicitly as:
\begin{equation}\label{eq_system}
  \begin{aligned}
    \frac{1}{4} \rho + \cK_{E}[\rho, \sigma, a_j] &= -\cS_0 \surfdiv
    \bEin, 
         &\quad &\text{on } \Gamma, \\
    \frac{1}{2} \sigma + \cK_{H} [\rho, \sigma, a_j] &= - \nhat \cdot \bHin,  
            & &\text{on } \Gamma, \\
     \int_{C_j} \bEtot \cdot d\bl &= 0, & &\text{for } j = 1,\ldots,
     2g,
  \end{aligned}
\end{equation}
with the operators~$\cK_E$ and~$\cK_H$ given by
\begin{equation}\label{eq_koperators}
  \begin{aligned}
    \cK_E[\rho,\sigma,a_j] &= ik\cS_0 \surfdiv \cS_k \bJ - \cC_k \rho
    - \cS_0 \surfdiv \nabla \times \cS_k \bK, \\
    \cK_H[\rho,\sigma,a_j] &= ik\nhat \cdot \cS_k \bK 
    - \cS'_k \sigma + \nhat \cdot \nabla \times \cS_k \bJ,
  \end{aligned}
\end{equation}
where
\begin{equation}
  \cS_k'f(\bx) = \frac{\partial}{\partial n} \int_\Gamma G(\bx,\bx')
  \, f(\bx) \, da(\bx'),
\end{equation}
with $\partial/\partial n$ denoting differentiation in the direction
normal to $\Gamma$ at $\bx$ and $G$ denoting the Green's function for
the Helmholtz equation:
\begin{equation}
  G(\bx,\bx') = \frac{e^{ik|\bx-\bx'|}}{4\pi | \bx-\bx' |}.
\end{equation}
The operator $\cC_k$ denotes the principal value part of the operator
$\cS_0 \surflap \cS_k$, given by:
\begin{equation}\label{eq:sls}
  \begin{aligned}
    \cC_k f(\bx) &=  \text{p.v. } \int_\Gamma
    \frac{1}{4\pi|\bx-\bx'|} \surfdiv \surfgrad \int_\Gamma G(\bx',\bx'')
    f(\bx'') \, da(\bx'') \, da(\bx') \\
    &=  -\text{p.v. } \int_\Gamma
    \surfgrad' \frac{1}{4\pi|\bx-\bx'|} \cdot
    \int_\Gamma \surfgrad' G(\bx',\bx'') f(\bx'') \, da(\bx'') \, da(\bx'),
  \end{aligned}
  \end{equation}
where~$\surfgrad'$ denotes the surface gradient in the $\bx'$
variable. See~\cite{oneil2017} for a computation of a similar
operator, $\cS_0 \surflap \cS_0$ using Cald\'eron projectors.

In the above expression, for simplicity reasons,
the dependence of~$\bJ$ and~$\bK$ on the variables $\rho$ and $\sigma$
and the coefficients $a_j$ has been suppressed.
All integral operators above are interpreted in their on-surface
principal-value sense. It is worth pointing out that the third terms
in $\cK_E$ and $\cK_H$ are in fact order-zero operators with respect
to the quantities $\bJ$ and $\bK$. They are of order minus-one with regard
to the variables $\rho$ and $\sigma$, as can be seen from the
construction of $\bJ$ and $\bK$.

\begin{remark}
The system of integral equations in~\eqref{eq_pecinteq1} is uniquely
invertible with the restriction that $\rho$, $\sigma$ (the generalized
Debye sources) and the data are all mean-zero functions. Without
explicitly enforcing this condition the system is rank-one deficient,
as can be seen from the operator $\cS_0 \surflap \cS_k$ (which has a
null-space of rank one). In
Section~\ref{sec_direct} we show how to incorporate this
constraint as a rank-one update to the system
matrix~\cite{sifuentes_2015}.
\end{remark}

\subsection{Dielectric bodies}
\label{sec_dielectric}

If the object $\Omega$ is not perfectly conducting,
but rather made of a dielectric material~\cite{jackson}, 
then the interfaces between
regions do not support physical current or charge sheets. The
materials merely become polarized, without the flow of charge,
resulting in conditions at interfaces based on continuity of the
normal and tangential components of the electric and magnetic
fields. These boundary conditions can be written as:
\begin{equation}
  \begin{aligned}
    \left[ \bn \times \bcE \right] &= 0, &\qquad \left[ \bn \cdot
      \epsilon \bcE \right] &= 0, \\
    \left[ \bn \times \bcH \right] &= 0, & \left[ \bn \cdot
      \mu \bcH \right] &= 0.
  \end{aligned}
\end{equation}
Recall that the quantities~$\bcD=\epsilon\bcE$ and~$\bcB=\mu\bcH$ are
the electric displacement and magnetic induction vector fields,
respectively; it is these fields whose normal components are
continuous across interfaces~\cite{papas}.

As in~\eqref{eq_gendebye}, we can use the generalized Debye source
representation to construct $\bE$, $\bH$ in each dielectric region
using layer potentials supported on the interface~\cite{EpGrOn}:
\begin{equation}
  \begin{aligned}
    \bE_\ell &= \sqrt{\epsilon_\ell \mu_\ell}
    \left( ik_\ell \bA_\ell - \nabla \varphi_\ell - \nabla \times
      \bQ_\ell \right), \\
    \bH_\ell &= \sqrt{\epsilon_\ell \mu_\ell}
    \left( ik_\ell \bQ_\ell - \nabla \psi_\ell + \nabla \times
      \bA_\ell \right).
  \end{aligned}
\end{equation}
Along the boundary of each dielectric region~$\Omega_\ell$ we define
tangential vector fields~$\bJ_\ell$, $\bK_\ell$ which will be
constructed to automatically satisfy the continuity conditions:
\begin{equation}
  \surfdivl \bJ_\ell = i k_\ell \rho_\ell,  \qquad
  \surfdivl \bK_\ell = i k_\ell \sigma_\ell.
\end{equation}
Let us denote by $\rho_{\ell'}$, $\sigma_{\ell'}$ the generalized
Debye sources supported on the \emph{opposite} side of
$\Gamma_\ell$ (i.e. those generalized Debye sources which are
responsible for generating fields in region~$\ell'$).
Then, as shown in~\cite{EpGrOn}, the following
construction of $\bJ_\ell$, $\bK_\ell$ leads to a uniquely invertible
system of second-kind integral equations. We describe the construction
in the case of only one component, $\Omega_1$, as it simplifies things
greatly. The case of multiple components (and nested components)
extends straightforwardly and is detailed in~\cite{EpGrOn}.
To this end, define the currents along~$\Gamma_0$ by:
\begin{equation}
  \begin{aligned}
    \bJ_0 &= i\omega \left( \sqrt{\epsilon_0 \mu_0}
      \surfgrad \surflap^{-1} \rho_0 - \bn \times \epsilon_1 \sqrt{\frac{\mu_1}
        {\epsilon_0}} \surfgrad\surflap^{-1} \rho_1 \right) + \bJ_0^h, \\
    \bK_0 &= i\omega \left( \sqrt{\epsilon_1 \mu_1}
      \surfgrad \surflap^{-1} \sigma_0 - \bn \times \mu_1
      \sqrt{\frac{\epsilon_1}
       {\mu_0}} \surfgrad\surflap^{-1} \sigma_1 \right) + \bK_0^h, 
  \end{aligned}
\end{equation}
with the (unrelated) harmonic components given by:
\begin{equation}
\bJ_0^h =  \sum_{j=1}^{2g} a_j \bh_j, \qquad
\bK_0^h =  \sum_{j=1}^{2g} b_j \bh_j,
\end{equation}
where $g$ denotes the genus of $\Gamma_0 = \Gamma_1$.  The Debye
currents $\bJ_1$, $\bK_1$ defining vector potentials in the interior
domain, $\Omega_1$, are then linearly related to the currents~$\bJ_0$
and~$\bK_0$, similar to the classic Muller integral
equation~\cite{muller}. Acting on tangential vector fields
along a boundary~$\Gamma$, given in their Hodge decomposition, 
let the linear operator~$\cU$ be defined as:
\begin{equation}
  \begin{aligned}
    \bF &= \surfgrad \alpha + \bn \times \surfgrad \beta + \bh, \\
    \cU\bF &= -\surfgrad \beta + \bn \times \surfgrad \alpha + \bh.
  \end{aligned}
\end{equation}
The operator~$\cU$ is known as
the \emph{clutching map}~\cite{EpGrOn}, and is equivalent to
$\cU = \cI$ when acting on harmonic components and $\cU = \bn \times$
when acting on the complement of the space of harmonic vector fields.
The interior Debye currents, $\bJ_1$, $\bK_1$, are then given as
\begin{equation}
\bJ_1 = \sqrt{\frac{\epsilon_0}{\epsilon_1}} \cU \bJ_0, \qquad
\bK_1 = \sqrt{\frac{\mu_0}{\mu_1}} \cU \bK_0.
\end{equation}
This construction ensures uniqueness
in the resulting integral equations.

It is important to note that by convention, we assume that
$\Gamma_\ell$ and $\Gamma_{\ell'}$ are oriented in the same manner,
i.e. that the unit normal vector along the boundary always points into
$\Omega_0$, the unbounded component of $\bbR^3$. Other conventions
would merely result in a change in some of the signs in the above
formulas. Also, note that our indexing of regions is opposite that
from~\cite{EpGrOn}. Here, we denote the outermost unbounded
region by $\Omega_0$, as is standard in the computational
electromagnetics literature.

Next, as shown in~\cite{EpGrOn}, the following boundary conditions
lead to a system of second-kind integral equations that are uniquely
invertible for all $\omega$, including in the limit as $\omega \to 0$:
\begin{equation}\label{eq_system2}
  \begin{aligned}
    \cS_0 \surfdiv \left[ \bEtot \right]  &= 0, & & \\
    \cS_0 \surfdiv \left[ \bHtot \right]  &= 0, & & \\
    \left[ \bn \cdot \sqrt{\epsilon} \bEtot \right]  &= 0, & & \\
    \left[ \bn \cdot \sqrt{\mu} \bHtot \right]  &= 0, &\quad 
          &\text{on } \Gamma, \\
     \int_{C_j} \left[ \bEtot \right] \cdot d\bl &= 0, & & \\
     \int_{C_j} \left[ \bHtot \right] \cdot d\bl &= 0, & 
         &\text{for } j = 1,\ldots,  2n_g,
  \end{aligned}
\end{equation}
where the contours~$C_j$, as before, span the space of the first
homology group in~$\bbR^3 \setminus \Omega$ and~$\Omega$.
The above system uniquely determines $\rho_\ell$, $\sigma_\ell$,
$\bJ_\ell^h$, and~$\bK_\ell^h$.
The system can be expanded as in Section~\ref{sec_pec}; we omit such
an expansion here for clarity as there are no new operators
introduced. However, it is worth pointing out that, in what follows, we
will always assume that the incoming electromagnetic fields~$\bEin$,
$\bHin$ are generated only from sources located
in~$\Omega_0$. In the case of the electric field and one
component, the interface conditions take the same form
as before:
\begin{equation}
  \begin{aligned}
     \cS_0 \surfdiv  \left( \bE_0 - \bE_1 \right) 
      &= -\cS_0 \surfdiv \bEin \\
    \bn \cdot \left( \sqrt{\epsilon_0} \bE_0 
      - \sqrt{\epsilon_1} \bE_1 \right) &= -\bn \cdot \sqrt{\epsilon_0} \bEin.
  \end{aligned}
\end{equation}
The jump conditions on the magnetic field are analogous.

In the next section,
we turn to a description of a separation of variables,
FFT-accelerated solver for the PEC system~\eqref{eq_system}
and dielectric system~\eqref{eq_system2} in the case of a single
component, $\Omega$, assuming it is a body of revolution.

\section{Integral equations on surfaces of revolution}
\label{sec_revolution}

The numerical solution of boundary integral equations on general
smooth surfaces in three dimensions has been an active area of
research for several decades, but robust, efficient and high-order accurate 
solvers are not easy to develop, mainly due to issues of surface 
representation and quadrature. However, when the geometry contains
rotational symmetry about, say, the $z$-axis, separation of variables
and Fourier decomposition can be used to turn the problem into a
sequence of boundary integral equations along a one-dimensional curve,
for which many fast and accurate numerical methods exist. 

\subsection{Geometry}
\label{sec_geometry}

In this section we establish the coordinate system to be used for a surface of
revolution~$\Gamma$, describe the corresponding Fourier analysis, and discuss
the form taken by surface differential operators. 
The scalar case has been carefully treated
in~\cite{young,helsing_2014}, which we extend here to the vector-valued case.
A special case of the present analysis, focusing on purely axisymmetric solutions
with applications to magnetohydrodynamics (MHD), 
was recently presented in~\cite{oneil2018taylor}. 

We denote the usual unit vector basis in Cartesian coordinates as
$\{ \ihat, \jhat, \khat \}$ and the usual unit vector basis in
cylindrical coordinates as $\{\rhat, \thetahat, \khat \}$:
\begin{equation}
  \begin{aligned}
    \rhat(\theta) &= \cos\theta \, \ihat + \sin\theta \, \jhat
    &\qquad \ihat &= \cos\theta \, \rhat(\theta) - \sin\theta \,
    \thetahat(\theta) \\
    \thetahat(\theta) &= -\sin\theta \, \ihat + \cos\theta \, \jhat
    &\qquad \jhat &= \sin\theta \, \rhat(\theta) + \cos\theta \,
    \thetahat(\theta)
    \end{aligned}
\end{equation}
where we have given formulae for both cylindrical and Cartesian unit
vectors in terms of the other. These relationships will be useful in
the following sections. The following \emph{addition formulae} will also be
useful:
\begin{equation}\label{eq_vecaddition}
  \begin{aligned}
    \rhat(\theta-\phi) &= \cos\phi \, \rhat(\theta) -
    \sin\phi \, \thetahat(\theta) \\
    \thetahat(\theta-\phi) &= \sin\phi \, \rhat(\theta) +
    \cos\phi \, \thetahat(\theta) \\
  \end{aligned}
\end{equation}
We will omit the implicit dependence (relative to $\ihat$ and $\jhat$)
of $\rhat$ and $\thetahat$ on the variable $\theta$ from now on,
unless it is explicitly required. In order to avoid confusion going
forward, {\em targets} will usually be denoted by $\bx$ and {\em
  sources} will usually be denoted by $\bx'$. This notation will be
consistent with integrating variables as well (integrals will be
performed in $r'$, $\theta'$, and $z'$ variables).

\begin{figure}[t]
  \centering
  \begin{subfigure}[b]{.35\linewidth}
    \centering
    \includegraphics[width=.95\linewidth]{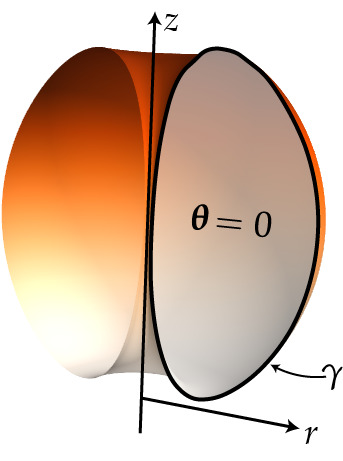}
    \caption{The generating curve~$\vct{\gamma}$.}
  \end{subfigure}
  \hfill
  \begin{subfigure}[b]{.6\linewidth}
    \centering
    \includegraphics[width=.95\linewidth]{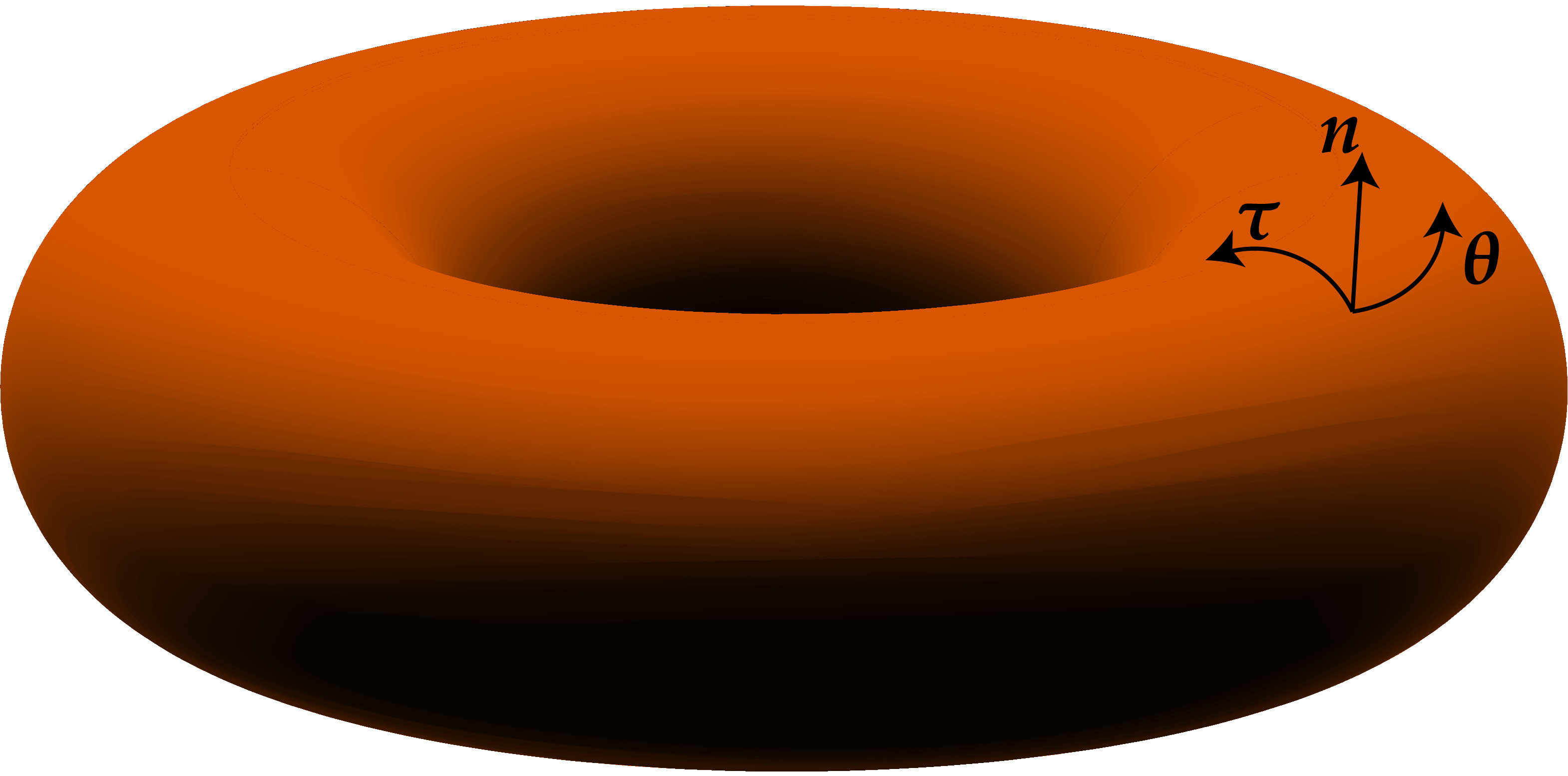}
    \caption{A local system of coordinates along~$\Gamma$.}
  \end{subfigure}
  \caption{Geometry of a body of revolution.}
  \label{fig_geo}
\end{figure}

Let the smooth closed curve $\vct{\gamma}: [0,L] \to \bbR^3$ be contained in
the $\theta=0$ plane, and let points on this curve be parameterized in
arclength as
\begin{equation}
  \vct{\gamma}(s) = \left( r(s), 0, z(s) \right),
\end{equation}
in cylindrical coordinates.  The boundary $\Gamma$
of a body of revolution $\Omega$ is then generated by rotating
$\vct{\gamma}$ about the $z$-axis.   A local
orthonormal basis along $\Gamma$ (assuming that $\vct{\gamma}$ is
oriented counter-clockwise) will be denoted by
$\{\tauhat, \thetahat, \nhat\}$ with $\thetahat$ the usual cylindrical
coordinate unit vector and
\begin{equation}
  \begin{aligned}
    \tauhat &= \frac{dr}{ds} \, \rhat + \frac{dz}{ds} \, \khat, \\
    \nhat &= \thetahat \times \tauhat .
  \end{aligned}
\end{equation}
See Figure~\ref{fig_geo} for a depiction of the setup and relevant
coordinate system. 

Lastly, we note that along a surface of revolution of genus 1
(i.e. a toroidal surface), 
the harmonic vector fields discussed in Section~\ref{sec_gendebye} are known
analytically and given by
\begin{equation} \label{hvecsformula}
\bh_1 = \frac{1}{r} \tauhat , \qquad \bh_2 = -\frac{1}{r} \thetahat.
\end{equation}

\subsection{The scalar function case}
\label{sec_scalar}

Scalar-valued boundary integral equations on surfaces of revolution
can be Fourier decomposed straightforwardly,
as in~\cite{young,helsing_2014}.  We
briefly state the decomposition here, and then discuss the
vector-valued case in detail.

Along the boundary $\Gamma$ of a body of revolution, a second-kind
integral equation with translation invariant 
kernel~$G$  takes the form:
\begin{equation}\label{eq_scalar}
\sigma(\bx) + \int_\Gamma G(\bx - \bx') \, \sigma(\bx') \, da(\bx')
 = f(\bx),
\end{equation}
for $\bx \in \Gamma$ and $da(\bx') = r' \, d\theta' \, ds'$.
Since~$G$ is assumed to be translation invariant,
this means that $G$, in particular, is rotationally
invariant: $G(\bx-\bx') = G(r,z,r',z',\theta-\theta')$.
Expanding the solution
$\sigma$ and right-hand side $f$ in terms of their Fourier series in
the azimuthal direction~$\theta$, we
have:
\begin{equation}\label{eq_sigfour}
\sigma(s,\theta) = \sum_{m=-\infty}^{\infty} \sigma_m(s) \,
e^{im\theta}, \qquad 
f(s,\theta) = \sum_{m=-\infty}^{\infty} f_m(s) \, e^{im\theta},
\end{equation}
where the Fourier coefficient functions $\sigma_m$ and $f_m$ are given
by
\begin{equation}
\sigma_m(s) = \frac{1}{2\pi} \int_{-\pi}^\pi \sigma(s,\theta) \,
e^{-im\theta} \, d\theta, \qquad 
f_m(s) = \frac{1}{2\pi} \int_{-\pi}^\pi f(s,\theta) \,
e^{-im\theta} \, d\theta.
\end{equation}
Substituting the Fourier representations for $\sigma$ and $f$
in~\eqref{eq_sigfour} into the integral equation in~\eqref{eq_scalar},
and collecting terms mode-by-mode, we have the following decoupled
integral equations along the generating curve $\vct{\gamma}$:
\begin{multline}
  \sigma_m(s) + \int_\gamma \sigma_m(s') \int_{-\pi}^\pi G(s,s',\theta-\theta')
  \, e^{-im\theta'} \, r' \, d\theta' \, ds' = f_m(s), \\
  \implies \quad \sigma_m(s) + 2\pi \int_\gamma G_m(s,s') \, \sigma_m(s') \,
  r' \, ds' = f_m(s)
\end{multline}
for all $m$, where we denote the Fourier modes of $G$ as $G_m$:
\begin{multline}
G(r,z,r',z',\theta-\theta') = \sum_{m=-\infty}^\infty G_m(r,z,r',z')
\, e^{im(\theta-\theta')} \\  \iff \quad G_m(r,z,r',z') =
  \frac{1}{2\pi}
  \int_{-\pi}^\pi G(r,z,r',z',\phi) \, e^{-im\phi} \, d\phi.
\end{multline}
We use the notation $G(s,s',\theta-\theta') = G(\bx-\bx')$  when
both $\bx$ and $\bx'$ are located on $\Gamma$, and likewise for the
Fourier modes $G_m$.
This separation of variables
procedure has reduced the problem of solving a boundary integral
equation on a surface embedded in three dimensions to that of solving
a sequence of boundary integral equations along the curve~$\gamma$,
which is embedded in two
dimensions. Given an efficient method for evaluating the modal Green's
functions $G_m$, performing the synthesis and analysis of $\sigma$ and
$f$ using the Fast Fourier Transform (FFT)~\cite{briggs_1995} yields a
very efficient direct solver for the original problem. 
The efficient
evaluation of $G_m$, for both the Laplace and Helmholtz problems, 
is a long-studied problem and discussed in Section~\ref{sec_greens}.

We now turn our attention to an analogous problem, that of Fourier
decomposing a vector-valued integral equation along a surface of
revolution.  While the generalized Debye system in~\eqref{eq_system}
is defined in terms of scalar-valued unknowns, its discretization
requires the application of layer potential operators to vector
quantities. Slightly more care is required than in the
scalar case.

\subsection{The vector-valued case}

As in the scalar case, 
a second-kind vector-valued boundary integral equation with
translation invariant kernel $G$ is of the form:
\begin{equation}\label{eq_vecinteq}
  \bu(\bx) + \int_\Gamma G(\bx,\bx') \, \bu(\bx') \, da(\bx') =
  \bsf(\bx).
\end{equation}
In the electromagnetic setting, both $\bu$ and $\bm f$ are
typically tangential vector fields, so we address that case only.  
Fully three-dimensional surface vector fields are encountered in
fluid dynamics and elasticity, and many of the following
calculations extend directly to those problems, albeit with components lying in the
direction normal to~$\Gamma$ as well.

Since $\Gamma$ is a surface of
revolution, the tangential vector fields in~\eqref{eq_vecinteq}
can be written in terms of
their Fourier series, component-wise, with respect to the local basis
$\tauhat$, $\thetahat$:
\begin{equation}\label{eq_vecfour}
  \begin{aligned}
    \bu(s,\theta) &= \sum_m \left( u_{\tau,m}(s) \, \tauhat
    + u_{\theta,m}(s) \, \thetahat \right) e^{im\theta} ,\\
    \bsf(s,\theta) &= \sum_m \left( f_{\tau,m}(s) \, \tauhat
    + f_{\theta,m}(s) \, \thetahat \right) e^{im\theta}.
  \end{aligned}
\end{equation}
Assuming that~$\bff$ is square-integrable along~$\Gamma$, and
that~\eqref{eq_vecinteq} is a Fredholm second-kind integral equation,
it can be shown that the solution~$\bu$ is also square
integrable~\cite{atkinson_1997}.  Therefore, the Fourier
representations in~\eqref{eq_vecfour} are complete. Since we have
that~$G$ is rotationally invariant, plugging the Fourier
expansion for $\bu$ into the integral term in~\eqref{eq_vecinteq} we
have:
\begin{multline}
 \int_\Gamma G(\bx-\bx') \, \bu(\bx') \, da(\bx')\\ = \int_\gamma
 \int_{-\pi}^\pi G(s,s',\theta-\theta') \sum_m \left( u_{\tau,m}(s') \,
 \tauhat(s',\theta') + u_{\theta,m}(s') \, \thetahat(\theta') \right)
 e^{im\theta'} \, r' d\theta' \, ds'.
\end{multline}
We now calculate the projection of this integral onto the
$m$th Fourier mode with respect to the unit vectors
$\rhat(\theta)$, $\thetahat(\theta)$.  The unit vectors
$\rhat(\theta')$ and $\thetahat(\theta')$ {\em cannot} be pulled
outside of the integrals, however it is a direct consequence of the
addition formulas in equation~\eqref{eq_vecaddition} that:
\begin{equation}
  \begin{aligned}
    \int_{-\pi}^\pi G(\cdot,\theta-\theta') \, \rhat(\theta') \,
    e^{im\theta'} \, d\theta' &= e^{im\theta}
    \int_{-\pi}^\pi G(\cdot,\phi) \, \rhat(\theta-\phi) \,
    e^{-im\phi} \, d\phi \\
    &= e^{im\theta} \, \rhat(\theta) \int_{-\pi}^\pi G(\cdot,\phi) \,
    \cos\phi \,  
    e^{-im\phi} \, d\phi  \\
    &\qquad -
    e^{im\theta} \, \thetahat(\theta) \int_{-\pi}^\pi G(\cdot,\phi) \,
    \sin\phi  \, 
    e^{-im\phi} \, d\phi  \\
    &= e^{im\theta} \left(  G_m^{\cos}(\cdot) \, \rhat(\theta) -
        G_m^{\sin}(\cdot) \, \thetahat(\theta) \right) 
  \end{aligned}
\end{equation}
where
\begin{equation}
  \begin{aligned}
  G_m^{\cos}(\cdot) &= \int_{-\pi}^\pi G(\cdot,\phi) \, \cos\phi
  \, e^{-im\phi} \, d\phi, \\
  G_m^{\sin}(\cdot) &= \int_{-\pi}^\pi G(\cdot,\phi) \, \sin\phi
  \, e^{-im\phi} \, d\phi.
\end{aligned}
\end{equation}
Likewise, we have that
\begin{equation}
\int_{-\pi}^\pi G(\cdot,\theta-\theta') \, \thetahat(\theta') \,
e^{im\theta'} \, d\theta' =
e^{im\theta} \left(  G_m^{\sin}(\cdot) \, \rhat(\theta) +
        G_m^{\cos}(\cdot) \, \thetahat(\theta) \right).
\end{equation}
Recalling that $\tauhat = \frac{dr}{ds} \rhat + \frac{dz}{ds} \khat$,
the integral equation in~\eqref{eq_vecinteq} can be re-written mode by
mode for all~$m$, and component-wise as:
\begin{equation}
  \begin{aligned}
    \rhat &:&\quad u_{\tau,m} \, \frac{dr}{ds}  + \int_\gamma \left(
    G^{\cos}_m(\cdot,s') \, \frac{dr}{ds'} \, u_{\tau,m} +
    G^{\sin}_m(\cdot,s') \, u_{\theta,m} \right) \, r' \, ds'
     &= \frac{dr}{ds} \, f_{\tau,m} ,\\
    \thetahat &:& u_{\theta,m}  + \int_\gamma \left(
    G^{\cos}_m(\cdot,s') \, u_{\theta,m} -
    G^{\sin}_m(\cdot,s') \, u_{\tau,m} \, \frac{dr}{ds'} \right) \,
    r' \, ds'
     &= f_{\theta,m} ,\\
    \khat &:& u_{\tau,m} \, \frac{dz}{ds} + \int_\gamma 
    G_m(\cdot,s') \, \frac{dz}{ds'} \, u_{\tau,m} \, 
     r' \, ds' \,
    &= \frac{dz}{ds} \, f_{\tau,m}.
  \end{aligned}
\end{equation}
The previous integral equations, while decoupled with respect to
Fourier mode $m$ and cylindrical coordinate vector components $\rhat$,
$\thetahat$, $\khat$, are {\em not} decoupled with respect to the
components of $\bu$, $u_{\tau,m}$ and $u_{\theta,m}$. This is due to the
action of the integral operator on the local coordinate system,
$\tauhat$, $\thetahat$, $\nhat$.  Depending on the exact boundary
conditions, linear combinations of the above integral equations may be
taken to form an even more coupled, second-kind system (since it is
obvious that one could solve for $u_{\tau,m}$ first, and then for
$u_{\theta,m}$). The main point here is to demonstrate that the
vector-value case is more involved than the scalar case, and that
there are additional kernels that need to be evaluated: $G_m^{\cos}$
and $G_m^{\sin}$.  Furthermore, the integral equations that we will
actually be discretizing involve the composition of layer potential
operators. This case can be handled analogously.

\subsection{Surface differentials}
\label{sec_differentials}

The integral equations resulting from the generalized Debye
representations make explicit use of the Hodge decomposition of
tangential vector fields, and therefore require the application of
surface differential operators.
In this section we briefly give formulae for these
operators along the axisymmetric boundary~$\Gamma$: 
the surface gradient $\nabla_\Gamma$, the surface
divergence $\nabla_\Gamma \cdot$, and the surface Laplacian
$\surflap$. Furthermore, the construction of surface currents
$\bJ$ and $\bK$ requires the application of $\surflap^{-1}$. We
provide the formula for $\surflap$ below, and then in
Section~\ref{sec_surflap} discuss a numerical method for inverting it.
In order to somewhat simplify the resulting
expressions, we assume that the generating curve $\gamma$ of the
boundary~$\Gamma$ is parameterized by arclength~$s$.  Similar formulae
with arbitrary parameterizations 
or  curvilinear coordinates can be found in most differential
geometry or mathematical physics texts, see for 
example~\cite{nedelec,frankel}.

Let \mbox{$f(s,\theta) = g(s)\, e^{im\theta}$} be a scalar
function defined on $\Gamma$ and
\[
\bF(s,\theta) = (F_\tau(s) \, \btau + F_\theta(s) \, \btheta) \,
  e^{im\theta}
\]
be a vector field tangential to $\Gamma$. We then have the following
three expressions for the surface differentials:
\begin{equation}\label{eq-surfops}
  \begin{gathered}
    \nabla_{\Gamma,m} f = \frac{\partial f}{\partial s} \, \btau  + \frac{im}{r}
f \, \btheta , \\
\nabla_{\Gamma,m} \cdot \bF = \frac{\partial F_\tau}{\partial s} +
\frac{1}{r} \frac{\partial r}{\partial s} F_\tau  + \frac{im}{r}
F_\theta, \\
\surflapm f = \frac{\partial^2 f}{\partial s^2} + 
\frac{1}{r} \frac{\partial r}{\partial s} \frac{\partial f}{\partial
  s} - \frac{m^2}{r^2} f.
\end{gathered}
\end{equation}
The on-surface quantity $\nabla \times \bF$, where $\nabla \times $
denotes the three-dimensional curl operator, can be computed for the
tangential $\bF$ above by using the standard curl in cylindrical
coordinates along with the relationships:
\begin{equation}
    \rhat = \frac{dr}{ds} \tauhat + 
    \frac{dz}{ds} \nhat, \qquad 
    \khat = \frac{dz}{ds} \tauhat -
    \frac{d r}{ds} \nhat. \\
\end{equation}

\section{A high-order solver}
\label{sec_solver}

In this section we detail the various parts of the solver,
including the efficient numerical evaluation of the modal Green's
functions, application of the surface differentials, inversion of
the surface Laplacian, and stable evaluation of the auxiliary
boundary conditions required by the Debye formulation.
The resulting linear systems, one for each Fourier mode, are then
solved directly using \textsc{LAPACK}'s dense LU routines.
This part of the computation is
not a dominant one, so no large effort was made to accelerate it.
First, we describe the Nystr\"om discretization of the integral equation.

\subsection{Nystr\"om discretization}
\label{sec_nystrom}

We restrict our examples in this work to integral equations along
smooth boundaries~$\Gamma$ of genus one.  As a consequence, the
generating curve $\gamma$, is represented by a smooth periodic, $\bbR^2$-valued function.
The resulting integral equations along~$\gamma$ are discretized at
equispaced points, and the resulting weakly-singular quadrature is
performed using hybrid Gauss-Trapezoidal
quadratures~\cite{alpert}. These quadrature rules are designed for
logarithmic singularities; while the fully three-dimensional Green's
function for the Helmholtz equation has a $1/r$ singularity, each of
the modal Green's functions only has a logarithmic singularity (as can
be seen from their special function representations~\cite{cohl_1999,conway_cohl}), 
and therefore high-order accuracy
can be obtained using these rules.  A thorough discussion of these
quadrature rules, and Nystr\"om discretization in general, is
contained in~\cite{hao_2014}, along with a comparison with other
methods.  Using this discretization method, a continuous integral
equation of the form
\begin{equation}
  \sigma(s) + \int_\gamma K(s,s') \, \sigma(s') \, ds' = f(s),
\end{equation}
is approximated by the finite dimensional system
\begin{equation}
  \sigma_j + \frac{2\pi L}{N} \sum_{k\neq j} w_{jk} \, K(s_j,s_k) \, \sigma_k=
  f(s_j), \qquad \text{for } j = 1,\ldots,N,
\end{equation}
where $L$ is the length of~$\gamma$, the $w_{jk}$ are the quadrature
weights for the Gauss-Trapezoidal rule, and $\sigma_j$ is an
approximation to the true solution $\sigma(s_j)$. The discretization
points are given as~$s_j = (j-1)L/N$.

In the case of an iterated integral, as in our case, each operator is
discretized separately and the resulting dense matrices are directly
multiplied together. For example, the continuous integral equation
\begin{equation}
  \sigma(s) + \int_\gamma K_2(s,s') \int_\gamma K_1(s',s'') \,
  \sigma(s'') \, ds'' \, ds' = f(s),
\end{equation}
is approximated by the finite-dimensional system
\begin{equation}
  \sigma_j + \frac{2\pi L}{N} \sum_{k\neq j} w_{jk} \, K_2(s_j,s_k) \, 
  \frac{2\pi L}{N} \sum_{k'\neq k} w_{kk'} \, K_1(s_k,s_{k'}) \, \sigma_k'=
  f(s_j), \qquad \text{for } j = 1,\ldots,N.
\end{equation}

Lastly, we note that in the numerical experiments contained in
Section~\ref{sec_examples} the integral equations for the generalized
Debye systems are only discretized using an \emph{odd} number of
points along the generating curve.
This avoids a spurious null-space issue arising
from the use of trapezoidal-based discretization schemes for Cauchy or
Hilbert-type integral operators, i.e. the null-space associated with
the alternating-point trapezoidal rule~\cite{sidi1988}.
Integral operators of this form arise
in the generalized Debye formulation due to the presence of surface
gradients and surface divergences. Furthermore, as mentioned above,
the discretization of iterated Hilbert-like integral operators, as
those in~\eqref{eq:sls}, requires repeated application of the
Gauss-Trapezoidal quadrature rule. Since this rule is a symmetric one,
each integral operator is accurately discretized, and therefore so is
the composition.
The analogous problem on fully 3D
surfaces requires extra care as well~\cite{bremer_2013}.

We now turn to the evaluation of the modal Green's functions.

\subsection{Modal Green's function evaluation}
\label{sec_greens}

The modal separation of variables calculation for integral equations
along surfaces is valid for any rotationally invariant Green's
function $G$.  Here we detail the precise calculations that are
necessary 
when $G$ is the Green's function for the Helmholtz equation:
\begin{equation}\label{eq_gk}
G(\bx,\bx') = \frac{e^{ik|\bx-\bx'|}}{4\pi |\bx-\bx'|}.
\end{equation}
Evaluation of the so-called {\em modal Green's functions} has received
varying degrees of attention over the past several decades, ranging
from brute-force integration to FFT-based
methods~\cite{gedney1990use,young}, as well as singularity
subtraction~\cite{helsing_2014}, exact special function
expansions~\cite{conway_cohl}, and contour
integration~\cite{gustafsson2010accurate}. While contour integration
for small modes $m$ is an attractive option due to the efficiency of
the resulting schemes, constructing deformations that yield accurate
results for higher modes is difficult, if not intractable.
Our evaluation scheme is based on the observations
in~\cite{helsing_2014}, with some additional numerical details filled in.

Switching to cylindrical coordinates and expanding in Fourier series,
\begin{equation}
  G(\bx,\bx') = \sum_{m=-\infty}^\infty G_m(r,z,r',z') \,
  e^{im(\theta-\theta')},
\end{equation}
where it is obvious that for $\bx \neq \bx'$:
\begin{equation}
G_m(r,z,r',z') = \frac{1}{2\pi} \int_{-\pi}^{\pi} G(r,z,r',z',\phi) \,
e^{-im\phi} \, d\phi.  
\end{equation}
Trivially, since $G$ is an even function, it is a function of
$\cos\phi$, and therefore
$G_{-m} = G_m$.
Furthermore, 
in fact, in the electro-/magneto-statics cases when $\omega=0$, the
Fourier modes can be calculated explicitly in terms of Legendre
functions of the second kind of half-order~\cite{cohl_1999,nist},
\begin{equation}\label{eq_qfuns}
  \frac{1}{4\pi| \bx - \bx |} = \frac{1}{4\pi^2 \sqrt{r r'}}
  \sum_m Q_{m-1/2}(\chi) \, e^{im(\theta-\theta')},
\end{equation}
where $\chi>1$ is given by
\begin{equation}
  \chi = \frac{r^2 + r'^2 + (z-z')^2 }{2rr'} .
\end{equation}
Here we are using notation consistent with previous discussions of
these functions \cite{cohl_1999,young}.
In the general case,~$\omega \neq 0$, we begin by 
expanding the integral representation of $G_m$ as
\begin{equation}\label{eq_gm}
\begin{aligned}
G_m(r,z,r',z') &= \frac{1}{2\pi} \int_{-\pi}^\pi 
\frac{e^{ik\sqrt{r^2 + r'^2 - 2rr' \cos\phi + (z-z')}}}
{4\pi \sqrt{r^2 + r'^2 - 2rr' \cos\phi + (z-z')}} \, e^{-im\phi} \,
d\phi \\
&= \frac{1}{8\pi^2 R_0} \int_{-\pi}^\pi 
\frac{e^{i\kappa\sqrt{1 - \alpha \cos\phi }}}
{ \sqrt{1 - \alpha \cos\phi }} \, e^{-im\phi} \,
d\phi
\end{aligned}
\end{equation}
where the parameters $\kappa$, $\alpha$, and $R_0$ are given by
\begin{equation}
\begin{gathered}
\kappa = k R_0, \qquad \alpha = 2rr'/R_0^2, \\
R_0^2 = r^2 + r'^2 + (z-z')^2,
\end{gathered}
\end{equation}
and we use $\phi = \theta-\theta'$ to denote the difference in
azimuthal angle.
In this form, it is obvious that any 
numerical scheme for evaluating $G_m$
will depend on three parameters: $\kappa$, $\alpha$, and $m$, with
$\alpha$ being the parameter which determines the growth of the
integrand near $\phi = 0$ 
($\alpha=1$ corresponds to the singularity in the Green's
function).
Based on several numerical experiments performed in quadruple
precision using the Intel Fortran compiler, and the methods developed
in~\cite{helsing_2014, young}, the following scheme provides near
machine accuracy evaluation of the sequence $G_{-M},\ldots,G_M$ in
roughly $\cO(M \log M)$ time using the FFT:


\begin{quote}
\textbf{For $\bm{0 \leq \alpha < 1/1.005 \approx .995}$}: These values
of $\alpha$ correspond to sources and targets that are reasonably
well-separated, and therefore the integrand in~\eqref{eq_gm} is
relatively smooth. For fixed $\kappa$, $\alpha$, and $R_0$, the
integral in~\eqref{eq_gm} can be discretized using the periodic
trapezoidal rule and
the kernels $G_{-M},\ldots,G_{M}$ and their
gradients can be computed using an $L$-point FFT, where $L$ depends on
$\kappa$ and $M$, i.e. the frequency content of the integrand.  If
$|\kappa| \leq 256$, a 1024-point FFT obtains near double precision
accuracy in this regime. If $|\kappa|>256$, an $L$-point FFT, with
$L\geq 2^{\log_2 |\kappa| + 2}$, obtains near double precision
accuracy.  Of course, $L$ must be chosen large enough, with $L > 2M$,
in order to compute $G_{-M},\ldots,G_{M}$. In practice, $M$ is usually much
smaller than the value of $L$ computed using the above estimates.
\end{quote}

\begin{quote}
\textbf{For $\bm{1/1.005 \leq \alpha < 1 }$}: In this case, the source
and target are relatively near to each other and the integrand
in~\eqref{eq_gm} is nearly singular. In this case, we follow the
method of~\cite{helsing_2014} and rewrite the integral as
\[
G_m(r,z,r',z') = \frac{1}{2\pi} \int_{-\pi}^\pi 
\frac{\cos(i\kappa \sqrt{1 - \alpha\cos\phi }) 
+ i \sin(i\kappa \sqrt{1 - \alpha\cos\phi })}
{4\pi R_0 \sqrt{1-\alpha\cos\phi}}
\, e^{-im\phi} \, d\phi.
\]
Letting
\[
H^c(\phi; \kappa,\alpha) = \frac{\cos(\kappa \sqrt{1 - \alpha\cos\phi }) }
{ \sqrt{1-\alpha\cos\phi}}, \qquad
H^s(\phi; \kappa,\alpha) = \frac{\sin(\kappa \sqrt{1 - \alpha\cos\phi }) }
{ \kappa \sqrt{1-\alpha\cos\phi}},
\]
we now note that $H^s$ is a very smooth function of $\phi$, even as
$\alpha \to 1$. Its Fourier modes can be computed using a
modestly~sized FFT, with size at most proportional to $\kappa$.

The Fourier modes of $H^c$ can be computed as the linear convolution
of the Fourier modes of $\cos(\kappa \sqrt{1 - \alpha\cos\phi }) $ and
the Fourier modes of $1/\sqrt{1-\alpha\cos\phi}$. The former modes can
be computed with the FFT, and the latter functions are known
analytically: they are proportional to $Q_{m-1/2}$, as
in~\eqref{eq_qfuns}. This method for evaluating $G_m$ was introduced
in~\cite{young}. See~\cite{briggs_1995} for more information on computing
linear convolutions with the FFT.
\end{quote}

We will not discuss the first case, when $0 \leq \alpha < 1/1.005$, as
the numerical evaluation is as simple as applying the FFT.
However,
there are some details which we will highlight concerning the more
difficult case, when $1/1.005 \leq \alpha < 1$.  In order to compute
$H^c$ via linear convolution of $Q_{m-1/2}$ with the Fourier modes of
$\cos(\kappa\sqrt{1-\alpha \cos\phi})$, one must first evaluate the
functions $Q_{m-1/2}(\chi)$.  The domain of relevance for these functions is
$(1,\infty)$, and they exhibit $\log$-type singularities at $\chi=1$.
These functions must be computed via backward-recurrence using
Miller's algorithm, as they are the recessive solution of the pair
$P_{m-1/2}$, $Q_{m-1/2}$ with respect to the index
$m$~\cite{gil_2007}.  Letting $\lambda = \sqrt{2/(\chi+1)}$, the first
two functions in this sequence are given explicitly in terms of
elliptic integrals~\cite{nist}:
\begin{equation}
 Q_{-1/2}(\chi) = \lambda \, K(\lambda), \qquad Q_{1/2}(\chi) = 
\chi \lambda  \, K(\lambda) 
 - \sqrt{2(\chi+1)} \, E(\lambda),
\end{equation}
where $K$ and $E$ are the complete elliptic integrals of the first-
and second-kind, respectively (in the notation of~\cite{nist}).
  Subsequent
terms obey the three-term recurrence
\begin{equation}
  \label{eq_gmrecurrence}
  Q_{m+1/2}(\chi) =  2\chi \frac{2m}{2m+1} Q_{m-1/2}(\chi) -
  \frac{2m-1}{2m+1} Q_{m-3/2}(\chi).
\end{equation}
However, for values of $\chi \approx 1$, near the singularity, an
excessive number of terms is required in order to run the forward
recurrence sufficiently far in Miller's algorithm.  For these values,
the forward recurrence is only mildly unstable and can be used with
care. Table~\ref{tab_fwdrec} provides estimates on the number of terms
that can be computed using the forward recurrence before values lose
more than 2 digits in absolute precision.  Using the fact that the
$Q_{m-1/2}$'s form a decreasing sequence, i.e.  $Q_{m+1/2} <
Q_{m-1/2}$, the values in Table~\ref{tab_fwdrec} were obtained
experimentally and set to be the index $m$ at which (numerically, in
double precision) $Q_{m+1/2} > Q_{m-1/2}$. This point in the sequence
indicates when loss of numerical precision begins. For values of $\chi >
1.005$, Miller's algorithm only requires only a few flops per
$Q_{m-1/2}$ for nearly full double-precision.

\begin{table}[!t]
  \begin{center}
    \caption{Sizes of forward recurrences that can be used for various
      values of $\chi$ to compute $Q_{m-1/2}(\chi)$ before
      significant loss of precision is encountered (more than
      $10^{-13}$ in absolute precision).}
  \label{tab_fwdrec}
    \begin{tabular}{|c|c|} \hline
      $\chi$  &  $\max m$ \\ \hline
      $(1.0,1.00000005]$ & 12307 \\ \hline
      $(1.00000005,1.0000005]$ & 4380 \\ \hline
      $(1.0000005,1.000005]$ & 1438 \\ \hline
      $(1.000005,1.00005]$ & 503 \\ \hline
      $(1.00005,1.0005]$ & 163 \\ \hline
    \end{tabular}
  \end{center}
\end{table}

Gradients of $G_m$ can be computed analytically using the integral
representation, and then evaluated using the previously discussed
methods. Using FFTs of the same size as when evaluating $G_m$ almost
always provides commensurate accuracy. Some sophistication can be used
when computing derivatives, as it involves products of several
terms. It is efficient to compute Fourier modes for each term and then
compute the overall convolution using the FFT. See~\cite{helsing_2014,
  young} for formulas regarding the gradients of $G_m$ and
$Q_{m-1/2}$.

Furthermore, the additional kernels required in the vector case,
$G^{\cos}_m$ and $G_m^{\sin}$, can be directly evaluated as
\begin{equation}
  \label{eq_modulated}
  \begin{aligned}
G_m^{\cos}  &= \int_{-\pi}^\pi G(\cdot,\phi) \, \cos\phi
\, e^{-im\phi} \, d\phi \\
&= \int_{-\pi}^\pi G(\cdot,\phi) \, \frac{e^{i\phi} + e^{-i\phi}}{2}
\, e^{-im\phi} \, d\phi \\
&= \frac{1}{2}\left( G_{m+1} 
+ G_{m-1} \right), \\
G_m^{\sin}  &= \frac{1}{2i}\left( G_{m+1} 
- G_{m-1} \right).
  \end{aligned}
\end{equation}
We now describe numerical tools related to surface differentials
along~$\Gamma$.

\subsection{Application of differentials}
\label{sec_surflap}

Using the generalized Debye source representation requires the
discretization of several surface-differential operators, namely
$\surfdiv$, $\surfgrad$, and $\surflap$.
Integration by parts along~$\Gamma$ allow for all composition of
differential and layer-potential operators to be constructed as
pseudo-differential operators of  (at most) order-zero
by direct differentiation of the Green's
function. For example, along $\Gamma$,
\begin{equation}
  \begin{aligned}
    \cS_0 \surfdiv \bJ (\bx) &= \int_\Gamma G_0(\bx,\bx') \, \surfdivp \bJ(\bx') \,
    da(\bx') \\
    &= - \int_\Gamma \surfgrad' G_0(\bx,\bx') \cdot \bJ(\bx') \, da(\bx'), \\
  \end{aligned}
\end{equation}
where $\surfgrad'$ denotes the surface gradient and
$\surfdivp$ denotes the surface divergence with respect
with respect to the
variable~$\bx'$.
Likewise, 
\begin{equation}
  \begin{aligned}
    \cS_0 \surflap \cS_k \rho (\bx) &= \cS_0 \surfdiv \surfgrad \cS_k
    \rho (\bx) \\
    &= -\int_\Gamma \surfgrad' G_0(\bx,\bx') \cdot \int_\Gamma
    \surfgrad' G_k(\bx',\bx'') \, \rho(\bx'') \, da(\bx'') \, da(\bx')
  \end{aligned}
\end{equation}
Using these identities, and the following construction of the operator
$\surfgrad \surflap^{-1}$, it is possible to build the discretized
system matrix {\em without} any numerical differentiation, despite the
presence of so many surface differential operators.

When using the generalized Debye source representation, as described
in Section~\ref{sec_pec}, to construct the surface vector fields
$\bJ$, $\bK$ from the sources $\rho$, $\sigma$ it is necessary to
apply the inverse of the Laplace-Beltrami operator $\surflap$. On a
general surface, this is rather complicated and a topic of ongoing
research~\cite{imbertgerard_2017,frittelli_2016,oneil2017}.
However, on a surface of revolution, the
procedure is reduced to solving a series of uncoupled ODEs with
periodic boundary conditions.
In order
to evaluate $\alpha = \surflap^{-1} f$ we instead solve the forward
problem
\begin{equation} \label{eq_integrodiff}
\begin{aligned}
\surflap \alpha &= f, \\
\int_\Gamma \alpha \, da &= 0,
\end{aligned}
\end{equation}
where it is assumed that $f$ is a mean-zero function.
On a surface of revolution, this can be reduced to a sequence of
uncoupled variable coefficient constrained ODEs, one for each mode
$m$:
\begin{equation}\label{eq_surfode}
\begin{aligned}
\left( \frac{\partial^2 }{\partial s^2} + 
\frac{1}{r} \frac{\partial r}{\partial s} \frac{\partial }{\partial
  s} - \frac{m^2}{r^2} \right)\alpha_m &= f_m, \\
\int_\gamma \alpha_m \, r \, ds &= 0,
\end{aligned}
\end{equation}
where, as before,
\begin{equation}
  \alpha(s,\theta) = \sum_m \alpha_m(s) \, e^{im\theta}, \qquad
  f(s,\theta) = \sum_m f_m(s) \, e^{im\theta},
\end{equation}
and~$r = r(s)$ denotes the $r$-coordinate along the generating
curve~$\gamma$.  As mentioned before, the surface Laplacian is
uniquely invertible as a map from mean-zero functions
to mean-zero functions \cite{EpGr}.  This mean-zero constraint is
automatically satisfied for right-hand sides $f_m$ with $m \neq 0$,
and the ODE portion of~\eqref{eq_surfode} is invertible.  For
rotationally symmetric functions~$f_0$, the integral constraint
in~\eqref{eq_surfode} must be explicitly enforced, as the ODE is not
invertible (trivially, the linear differential operator has a
null-space of constant functions).  One option for solving this system
is to discretize~$\surflap$ pseudo-spectrally as
in~\cite{imbertgerard_2017} and invert the resulting dense
matrix. This procedure may suffer from ill-conditioning when~$\gamma$
requires many discretization points.  Alternatively, letting
$u_m = d^2 \alpha_m /d s^2$, this ODE can formally be rewritten as a
second-kind integral equation with periodic boundary conditions:
\begin{equation}
\begin{gathered}
u_m + \frac{1}{r} \frac{\partial r}{\partial s} \int u_m \, ds -
\frac{m^2}{r^2} \iint  u_m \, ds' \, ds = f_m, \\
u_m(0) = u_m(L), \qquad \frac{d u_m}{d s}(0) = 
\frac{d u_m}{d s}(L).
\end{gathered}
\end{equation}
Here, we denote by~$\int u \, ds$ the (mean-zero) anti-derivative of~$u$.  The
previous integral equation can be solved for $u_m$ by using the
discrete Fourier transform (DFT)
and enforcing the proper integral constraint explicitly on
$u_0$. See~\cite{trefethen_spectral_book} or~\cite{wright2015} for a
discussion of methods for solving ODEs with periodic boundary
conditions.

Once $u_m$ (or its Fourier series) has been computed, $\alpha_m$ and
its first derivative can be easily and stably computed in the Fourier
domain by division of the mode number.
Both~$\alpha_m$ and its first derivative are needed to compute
its surface gradient, as per the surface differential
formulas in~\eqref{eq-surfops}.


\subsection{Auxiliary conditions for surfaces of nontrivial genus}
\label{sec_auxiliary}

The generalized Debye formulation for scattering from genus 1 (or
higher) objects involves a finite set of
integral constraints on field quantities in order to fix the
projection of $\bJ$, $\bK$ onto harmonic vector fields, as per the
Hodge decomposition in~\eqref{eq_hodge}. In the case of scattering
from a genus 1 PEC, we enforce the additional constraints 
\begin{equation}\label{eq_cycles}
\int_{C_A} \bEtot \cdot d\bl  = 0, \qquad 
\int_{C_B} \bEtot \cdot d\bl = 0,
\end{equation}
where the loops $C_A$ and $C_B$ are shown in Figure~\ref{fig_cycles}.
These integral conditions come from the fact that along the surface of
a PEC, it must be the case that $\bn \times \bcE = \vct{0}$ pointwise; in
particular, this means that any line integral of the tangential
components of the electric field must be zero.
In the case of scattering from a penetrable object of genus 1, four
additional constraints must be enforced:
\begin{equation}
  \begin{aligned}
    \int_{C_A} \left[ \bEtot \right] \cdot d\bl  &= 0, &\qquad 
    \int_{C_B} \left[ \bEtot \right] \cdot d\bl  &= 0, \\
    \int_{C_A} \left[ \bHtot \right] \cdot d\bl  &= 0, &\qquad 
    \int_{C_B} \left[ \bHtot \right] \cdot d\bl  &= 0.
  \end{aligned}
\end{equation}
We provide details for discretizing the integrals
in~\eqref{eq_cycles}, as the integrals in the dielectric case are
identical.  Along general complex surfaces embedded in 3D, the
computational task of finding \emph{loops} must first be performed,
and then the loops must be discretized. However, given that our
scattering geometry is a body of revolution, in the local coordinate
system $\tauhat, \thetahat, \nhat$ along $\Gamma$ these integrals can
be rewritten as:
\begin{equation}\label{eq_cycles2}
  \begin{aligned}
    \int_{C_A} \bE \cdot \tauhat \, ds &= - \int_{C_A} \bEin \cdot
    \tauhat \, ds, \\
    \int_{C_B} \bE \cdot \thetahat \, d\theta &= - \int_{C_B} \bEin \cdot
    \thetahat \, d\theta.
  \end{aligned}
\end{equation}

Trivially, if $\bEtot$ were decomposed into Fourier modes in the
azimuthal direction, then only the $m=0$ mode would contribute to the
line integrals along $C_B$. However, this is \emph{not} the case for
integrals along~$C_A$. In this case, it is likely that all Fourier
modes of~$\bEtot$ (in the variable~$\theta$)
contribute to the line integral. This is not
convenient for a separation of variables solver because it
effectively couples the modal integral equations (due to the fact that
the harmonic vector fields along~$\Gamma$ only generate purely
axisymmetric fields). These auxiliary conditions are needed to prove
uniqueness results for the systems of Fredholm equations that we solve 
below. The uniqueness argument continues to apply if
the line integral conditions on~$C_A$
in~\eqref{eq_cycles} and~\eqref{eq_cycles2} are replaced by a surface
integral that effectively integrates out the contribution from all 
non-zero Fourier modes:
\begin{equation}
\int_\Gamma \bEtot \cdot \btau \, da  = 0.
\end{equation}
In what follows this is the condition that we will use.

Consider now the integrals over a $B$-cycle.
As $k \to 0$, if the incoming field~$\bEin$,~$\bHin$
is due to sources exterior to a ball containing~$\Omega$, then the circulation of $\bEin$
along~$C_B$ is $\cO(k) \to 0$. Indeed, by
application of Maxwell's equations and Stokes's Theorem,
\begin{equation}
\int_{C_B} \bEin \cdot \thetahat \, ds = -\int_{S_B} ik \, \bHin
\cdot d\ba,
\end{equation}
where~$S_B$ is a spanning surface (in the exterior of~$\Omega$)
with boundary~$C_B$. 
As discussed in~\cite{epstein-2013,EpGrOn}, the same is true for the scattered field:
\begin{equation}\label{eq_scaling}
    \int_{C_B} \bE \cdot \thetahat \, d\theta \sim \cO(k).
\end{equation}
As a result of this scaling,
difficulties arise in accurately discretizing the circulation of~$\bE$
along~$C_B$ as both the integral and the data tend to zero.
Without care, catastrophic cancellation can result in a loss of
accuracy due to ill-conditioning of the system matrix.

In the case of the MFIE, this ill-conditioning can be avoided by 
adding a constraint on the
vector potential $\bA$ itself (see~\cite{epstein-2013}).
In the generalized Debye case, however, we are using non-physical variables
and need to stabilize the problem using $\bE$ and/or $\bH$ themselves.
For this, following the discussion of~\cite{epstein-2013}, we denote by~$\bE(k)$ the
electric field generated by the Debye sources~$\rho$, $\sigma$ and
currents~$\bJ$, $\bK$ with wavenumber~$k$:
\begin{equation}
  \bE(k) = ik \cS_k\bJ - \nabla  \cS_k \rho - \nabla \times
  \cS_k\bK.
\end{equation}
Furthermore, from~\eqref{eq_scaling}, we see that
\begin{equation}
  \lim_{k \to 0} \int_{C_B} \bE(k) \cdot \thetahat \, d\theta =
   \int_{C_B} \bE(0) \cdot \thetahat \, d\theta = 0.
\end{equation}
Using this fact, note that
for any~$k>0$,
\begin{equation}
  \begin{aligned}
    \frac{1}{k} \int_{C_B} \bE(k) \cdot \thetahat \, d\theta &=
    \frac{1}{k} \left( \int_{C_B} \bE(k) \cdot \thetahat \, d\theta  -
      \int_{C_B} \bE(0) \cdot \thetahat \, d\theta \right) \\
    &=
     \int_{C_B} \left(\frac{ \bE(k)-\bE(0)}{k} \right) \cdot \thetahat
     \, d\theta \\
     &=
     \int_{C_B} \left(\frac{ \bE(k)-\bE(0)}{k} \right) \cdot \thetahat
     \, d\theta \\
     &=
     \int_{C_B} \bE_{\text{diff}} \cdot \thetahat \, d\theta,
  \end{aligned}
\end{equation}
where
\begin{equation}
  \bE_{\text{diff}} = i\cS_k\bJ - i \nabla \times \cS_k  \left( \surfgrad
    \surflap^{-1} \rho - \nhat \times \surfgrad \surflap^{-1} \sigma
  \right)
  - \nabla \times \frac{1}{k} \left( \cS_k - \cS_0 \right)\bK_H.
\end{equation}
While catastrophic cancellation (in relative precision)
occurs in computing the $B$-cycle circulation of $\bE$, the same
circulation of $\bE_{\text{diff}}$ is $\cO(1)$ and can be computed
stably as long as care is taken in evaluation of the difference operator:
\begin{equation}
  \cS_{\text{diff}} = \frac{1}{k} \left( \cS_k - \cS_0 \right).
\end{equation}
For small values of~$k|\bx-\bx'|$,
expanding the numerator of the Green's function in a
Taylor series and taking the first several terms achieve
near machine precision. With this computation, instead of enforcing
the circulation condition $\int_{C_B} \bEtot$ we instead enforce the
identical condition:
\begin{equation}
  \begin{aligned}
    \int_{C_B} \frac{\bE(k) - \bE(0) }{k} \cdot d\bl &=
    -\frac{1}{k} \int_{C_B} \bEin \cdot d\bl \\
    &= \frac{1}{k} \int_{S_B}\nabla \times \bEin \cdot d\ba \\
    &= -i \int_{S_B} \bHin \cdot d\ba,
  \end{aligned}
\end{equation}
where $S_B$ is the spanning surface for the cycle $C_B$. In the case
of a body of revolution, the surface~$S_B$ is easily discretized using
the trapezoidal rule in~$\theta$ and Gauss-Legendre quadrature
in~$r$. This provides as spectrally accurate scheme for computing the
right-hand side as the incoming data is assumed to be smooth.

Next, we describe the full multi-mode direct solver for computing the
full electromagnetic scattering problem.

\subsection{A direct solver}
\label{sec_direct}

Using the above numerical machinery, a full system matrix discretizing
each modal integral equation can be constructed.
We describe this process in detail for the PEC case. The dielectric
problem is analogous but somewhat more involved.
Let $N$ denote the number of discretization points along the
generating curve $\gamma$ and $M/2$ denote the maximum Fourier mode
required in the azimuthal direction. The surface $\Gamma$ is then
sampled using an $N \times M$ tensor-product grid in $[0,L] \times
[0,2\pi]$. The fully 3D incoming fields
\begin{equation}
  \begin{aligned}
    \bEin &= \Ein_x \ihat + \Ein_y \jhat + \Ein_z \khat \\
    \bHin &= \Hin_x \ihat + \Hin_y \jhat + \Hin_z \khat
g  \end{aligned}
\end{equation}
are sampled on this grid and converted to a cylindrical coordinate
vector field:
\begin{equation}
  \begin{aligned}
    \bEin &= \Ein_r \rhat + \Ein_\theta \thetahat + \Ein_z \khat \\
    \bHin &= \Hin_r \rhat + \Hin_\theta \thetahat + \Hin_z \khat.
  \end{aligned}
\end{equation}
The Fourier series of each of these cylindrical components is computed
via an FFT so that:
\begin{equation}
  \Ein_r = \sum_{m = -M/2}^{M/2} \Ein_{r,m} \, e^{im\theta},\ 
  \Ein_\theta = \sum_{m = -M/2}^{M/2} \Ein_{\theta,m} \, e^{im\theta},\ 
  \Ein_z = \sum_{m = -M/2}^{M/2} \Ein_{z,m} \, e^{im\theta},\ {\rm etc.}
\end{equation}
The mode-by-mode tangential projections in the local orthonormal basis,
following the notation in \eqref{eq_vecfour},
are then computed as:
\begin{equation}
  \begin{aligned}
    \Ein_{\tau,m} &= \frac{dr}{ds} \Ein_{r,m} + \frac{dz}{ds}
    \Ein_{z,m}, &\qquad \Hin_{\tau,m} &= \frac{dr}{ds} \Hin_{r,m} + \frac{dz}{ds}
    \Hin_{z,m}, \\
    \Ein_{\theta,m} &= \frac{dz}{ds} \Ein_{r,m} - \frac{dr}{ds}
    \Ein_{z,m}, &\qquad \Hin_{\theta,m} &= \frac{dz}{ds} \Hin_{r,m} - \frac{dr}{ds}
    \Hin_{z,m}.
  \end{aligned}
\end{equation}
For each mode $m = -M/2,\ldots, M/2$, the projection of the data
$\cS_0 \surfdiv \bEin$ onto the $m$th Fourier mode in $\theta$
is computed as:
\begin{equation}
  \begin{aligned}
    \left( \cS_0 \surfdiv \bEin \right)_m &=
    \left(  \int_\Gamma G^0(\bx,\bx') \, \surfdiv \bEin(\bx') \, da(\bx')
    \right)_m \\
    &= \left(  -\int_\Gamma \surfgrad' G^0(\bx,\bx') \cdot
      \bEin(\bx') \, da(\bx') \right)_m \\
    &= -2\pi \int_\gamma \left( \frac{dG^0_m}{ds'}(s,s') \Ein_{\tau,m}(s') -
      \frac{im}{r'} G^0_m(s,s') \, \Ein_{\theta,m} \right) r' \, ds',
  \end{aligned}
\end{equation}
where for $\bx$ and $\bx'$ on $\Gamma$,
\begin{equation}
  G^0(\bx,\bx') = \frac{1}{4\pi |\bx-\bx'|} \qquad\text{and} \qquad
  G^0(\bx,\bx') = \sum_m G^0_m(s,s') \, e^{im(\theta-\theta')}.
\end{equation}
The $m$th Fourier mode of the scalar data $\bn \cdot \bHin$ is denoted by
$\Hin_{n,m}$. A sequence of decoupled integral equations on $\gamma$
can then be solved for each mode.
For this, we expand the Debye sources as
\begin{equation}\label{sourceexp}
\rho(s,\theta) = \sum_{m=-\infty}^{\infty} \rho_m(s) \,
e^{im\theta}, \qquad 
\sigma(s,\theta) = \sum_{m=-\infty}^{\infty} \sigma_m(s) \,
e^{im\theta}.
\end{equation}
For $m \neq 0$, we define $\bJ_m,\bK_m$ by 
\begin{eqnarray*}
\bJ_m(s) &=& ik \left(  \nabla_{\Gamma,m}  \surflapm^{-1} \rho_m -
 \nhat \times \nabla_{\Gamma,m}  \surflapm^{-1} \sigma_m \right) \\
\bK_m(s) &=& \nhat \times \bJ_m(s).
\end{eqnarray*}
For $\bJ_0,\bK_0$, we must include the contributions of the harmonic vector
fields:
\begin{eqnarray*}
\bJ_0(s) &=& ik \left(  \nabla_{\Gamma,0}  \Delta_{\Gamma,0}^{-1} \rho_0 -
 \nhat \times \nabla_{\Gamma,0}  \Delta_{\Gamma,0}^{-1} \sigma_0 \right) +
a_1 \bh_1 + a_2 \bh_2, \\
\bK_0(s) &=& \nhat \times \bJ_0(s),
\end{eqnarray*}
where the basis vectors $\bh_1,\bh_2$ are given by
\eqref{hvecsformula}.
For $m \neq 0$, we solve the Nystr\"om discretization of the
system
\begin{equation}\label{eq_modeeq}
  \begin{gathered}
 \frac{1}{4} \rho_m + ik\cS_{0,m} \nabla_{\Gamma,m}\cdot \cS_{k,m}
  \bJ_m - \cS_{0,m} \Delta_{\Gamma,m} \cS_{k,m} \rho_m +
  \cS_{0,m} \nabla_{\Gamma,m}\cdot \nabla
  \times \cS_{k,m} \bK_m =  \left( \cS_0 \surfdiv \bEin
\right)_m \\
\frac{1}{2} \sigma_m + ik \bn \cdot \cS_{k,m} \bK_m - \cS'_{k,m} \sigma_m +
\bn \cdot \nabla \times \cS_{k,m} \bJ_m = \Hin_{n,m},
\end{gathered}
\end{equation}
for the unknowns $\rho_m$ and $\sigma_m$.
The differential operators here are given 
in~\eqref{eq-surfops} and the integral operators $\cS_{k,m}$ are defined by
\begin{equation}
  \cS_{k,m}f(s) = 2\pi \int_\gamma G_m(s,s') \, f(s') \, r' \, ds',
\end{equation}
with~$G_m$ the Fourier modes of~$G$ in~\eqref{eq_gk}. It is implied
that each of the layer potential 
operators above has been discretized using
an~$N$-point Nystr\"om scheme with hybrid Gauss-trapezoidal quadrature
corrections~\cite{alpert}.

For the $m=0$ case, we have two additional unknowns (the coefficients $a_1,a_2$ of the
harmonic vector fields), and we augment the system of equations \eqref{eq_modeeq}, 
following \cite{EpGrOn}, with the constraints
\begin{equation}\label{eq_modeeq0}
  \begin{gathered}
\int_\Gamma E_{\tau,0} \, da = - \int_\Gamma \bEin \cdot \btau \, da, \\
2\pi r_B  \frac{E_{\theta,m}(k) - E_{\theta,m}(0) }{k}  = 
    -i \int_{S_B} \bHin \cdot d\ba.
\end{gathered}
\end{equation}
Here,~$r_B$ denotes the radius of the $B$-cycle, $C_B$. The
functions~$E_{\theta,m}(k)$ ,$E_{\theta,m}(0)$ are assumed to be
evaluated at the point along the generating curve which corresponds
to~$C_B$.

\section{Numerical examples}
\label{sec_examples}

In this section, we provide several numerical examples to demonstrate
the efficiency and accuracy of the electromagnetic scattering solvers
described in this work. The code was implemented using a mixture of Fortran
77/90, compiled with the Intel Fortran Compiler, and linked against
the Intel MKL libraries for low-level linear algebra routines. 
Minimal efforts were
made to parallelize the system matrix assembly using OpenMP
directives, but no fine-grained optimizations were made. Unless
otherwise noted, examples were run on a workstation with 32
Intel Xeon Gold 6130 cores at 2.1GHz with 512GB of shared memory.
Plots of 3D images were created in Paraview~\cite{paraview2005}.

In each of the examples, an estimate of the \emph{resolution of the
geometry} is provided. This is an estimate of how accurately the
parameterization of the geometry and its derivatives have been approximated
in the discretization. If $r_1, \ldots, r_n$ denotes a sample of the
$r$-component of the geometry (in arclength), then the resolution of
the $r$-component is estimated as follows. First, we compute the
discrete Fourier transform of the sequence:
\begin{equation}
  \hat{r}_k = \frac{1}{n} \sum_{j=-n/2}^{n/2} r_j \, e^{-2\pi i j k /n },
  \quad \text{for } k=-n/2, \ldots,n/2.
\end{equation}
Letting $\Vert \hat{r} \Vert$ denote the $\ell^2$ norm of $\hat{r}$,
the resolution of the sampling $r_1, \ldots, r_n$ is estimated as
\begin{equation}
  \Res{r} = \frac{\sqrt{\sum_{-n/2}^{-n/2+3} |\hat{r}_k|^2 +
\sum_{n/2-3}^{n/2} |\hat{r}_k|^2  }}{\Vert \hat{r} \Vert}.
\end{equation}
This estimate, $\Res r$, computes the energy located in the tail of
the Fourier series of $r$. Examining the last four positive and
negative terms is a robust heuristic, and other estimates could certainly be used. 
For each geometry, the above estimate is computed for each component of
$\gamma$ and $\gamma'$, where as before, $\gamma$ denotes the
generating curve of the full surface $\Gamma$. The resolution of
$\Gamma$, $\Res \Gamma$, is taken to be the maximum value of these
four individual resolutions.

%

\subsection{Convergence results}

In our first set of examples, the accuracy of the solver is
demonstrated by testing against an exact
solution. Exact solutions to Maxwell's equations 
in the exterior (or interior) of the object~$\Omega$ can be generated by
placing a small current loop in the interior (or exterior). 
By solving the Debye source integral equation with the corresponding boundary data, 
we can evaluate the field at an arbitrary point and compare with the known, exact
solution. For the exact solution, we define the
field due to a counter-clockwise oriented
current density supported on a horizontal loop centered
at~$\vct{c}_0$ with radius~$r_0$ as:
\begin{equation}
  \begin{aligned}
  \bH^{ex}(\bx) &= \nabla \times  \int_0^{2\pi} 
  \frac{e^{ik|\bx - \vct{c}_0 - r_0 \rhat(\theta)|}}
{|\bx - \vct{c}_0 - r_0 \rhat(\theta)|} \, 
\frac{\btheta(\theta)}{2\pi r_0} \, d\theta, \\
  \bE^{ex}(\bx) &= \frac{i}{k} \nabla \times \bH^{ex}(\bx),
\end{aligned}
\end{equation}
with $\rhat$, $\btheta$ denoting the unit vectors in 
cylindrical coordinates.
The above fields can easily be evaluated using analytic
differentiation combined with  a trapezoidal
discretization of the resulting integrals.

\begin{figure}[!t]
  \centering
  \begin{subfigure}[b]{.45\linewidth}
    \centering
    \includegraphics[width=.95\linewidth]{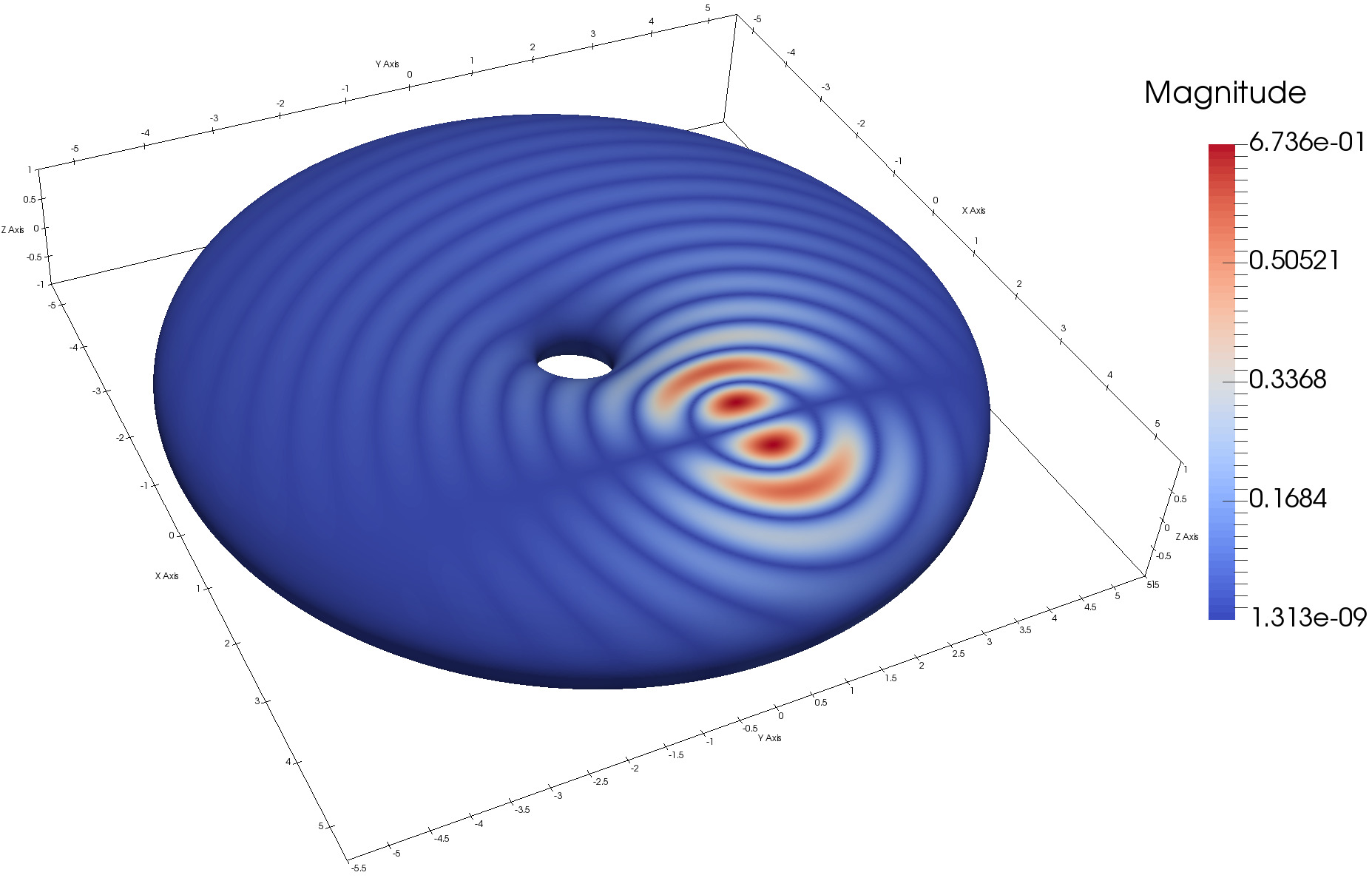}
    \caption{Test geometry and exterior data with $\epsilon_0 =
      39.48$ and $\lambda_0=1.0$.}
  \end{subfigure}
  \hfill
  \begin{subfigure}[b]{.45\linewidth}
    \centering
    \includegraphics[width=.95\linewidth]{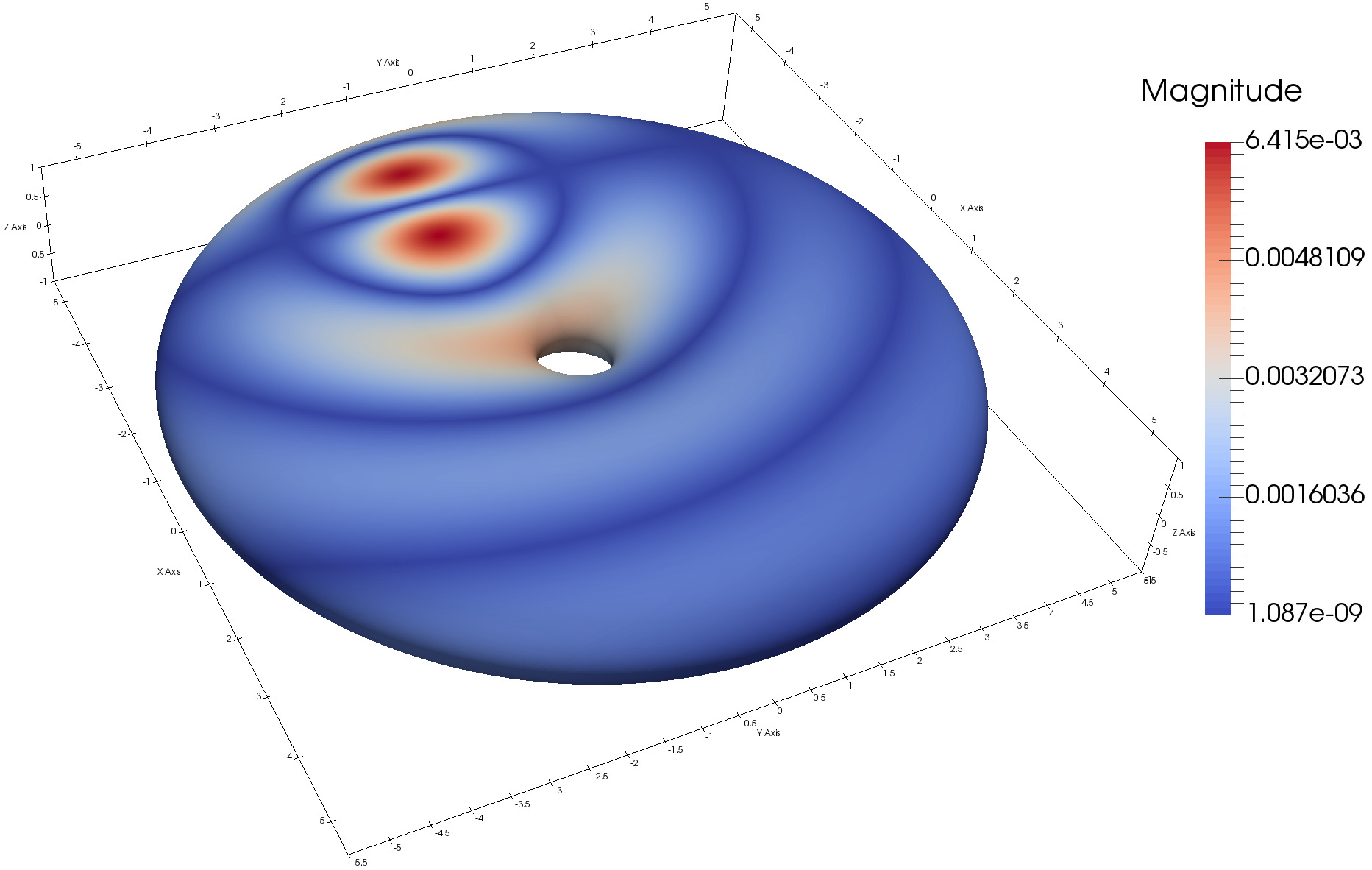}
    \caption{Test geometry and interior data with $\epsilon_1 =
      2.47$ and $\lambda_1=4.0$.}
  \end{subfigure}
  \caption{The geometry for dielectric convergence tests, a 2.5-1 torus
    centered at the origin. Shown is $|\realpart \bHin|$ along the surface
    of the scatterer for both the interior data and exterior data. In
    each case, $\omega = 1.0$ and $\mu=1.0$.}
  \label{fig_testgeo2}
\end{figure}

In the following
tables reporting convergence results, $N$ denotes the number of
discretization points used along the generating curve, $L$ denotes the
number of points in the azimuthal direction, $M = NL$, $t$ denotes the solution
time in seconds (not including the evaluation of the incoming field along
the surface),
$t_{G_m}$ denotes the time (in seconds) required to evaluate all the
modal Green's functions,
and
the error reported is the relative $\ell^2$ error in the computed field and the
known field on a sphere, as
estimated via 50 points equispaced in $\theta$ and $\phi$:
\begin{equation}
  Err(\bE) = \frac{\sum_i \| \bE(\bx_i) - \bE^{ex}(\bx_i)
    \|^2 }
  {\sum_i \|  \bE^{ex}(\bx_i) \|^2 }.
\end{equation}
Lastly, $\lambda$ denotes the
wavelength of the driving field, i.e. $\lambda = 2\pi/k$.
The sample geometries for the perfect electric conducting case
are shown in Figure~\ref{fig_testgeo}. The test geometries and data
for verifying the convergence of the dielectric problem are
shown in Figure~\ref{fig_testgeo2}. 
Convergence results are reported in Tables~\ref{tab_conv_pec}
and~\ref{tab_conv_die}.

\afterpage{
\begin{figure}[!t]
  \centering
  \begin{subfigure}[b]{.45\linewidth}
    \centering
    \includegraphics[width=.95\linewidth]{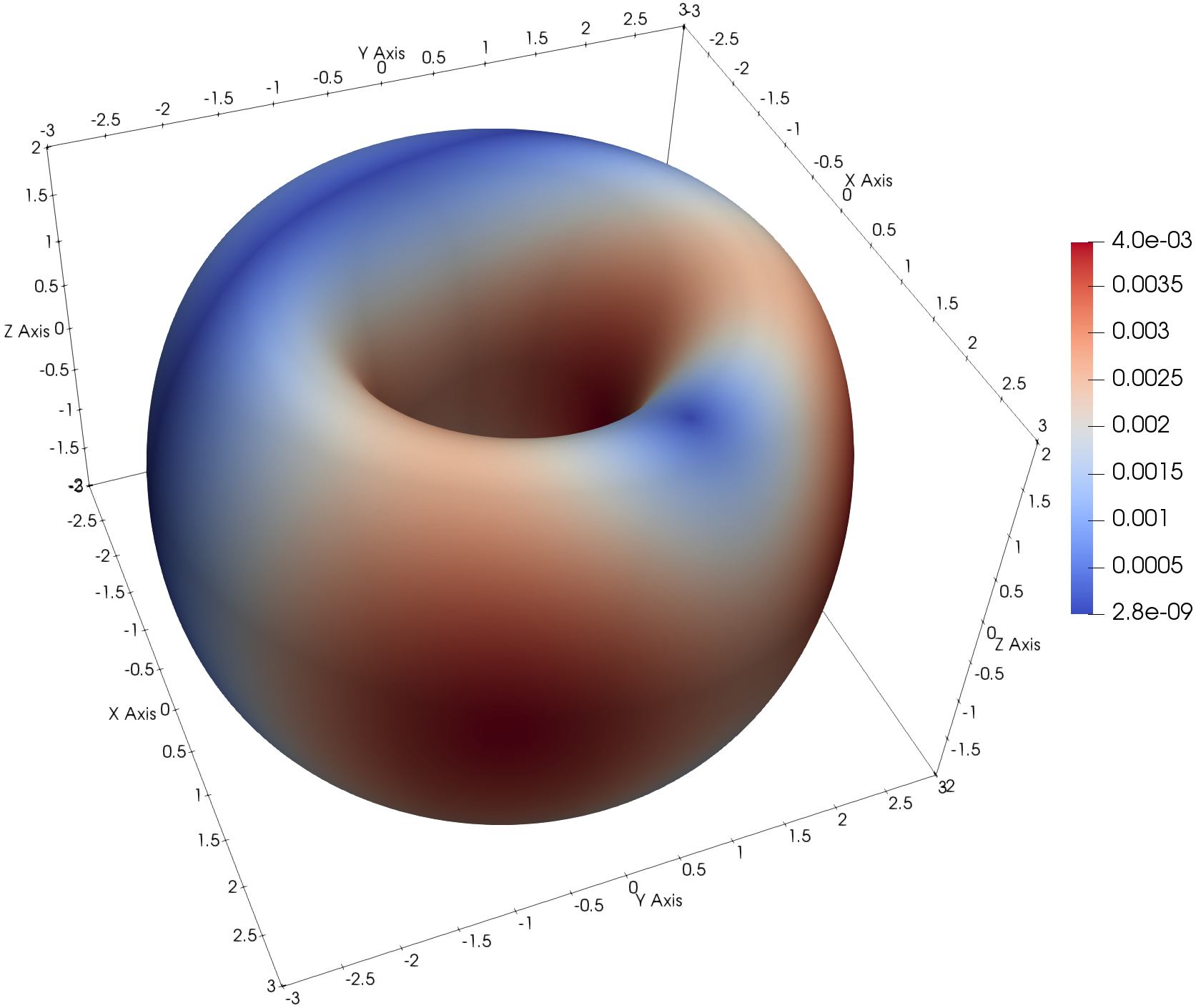}
    \caption{The test geometry and data for $\lambda=6.0$.}
  \end{subfigure}
  \hfill
  \begin{subfigure}[b]{.45\linewidth}
    \centering
    \includegraphics[width=.95\linewidth]{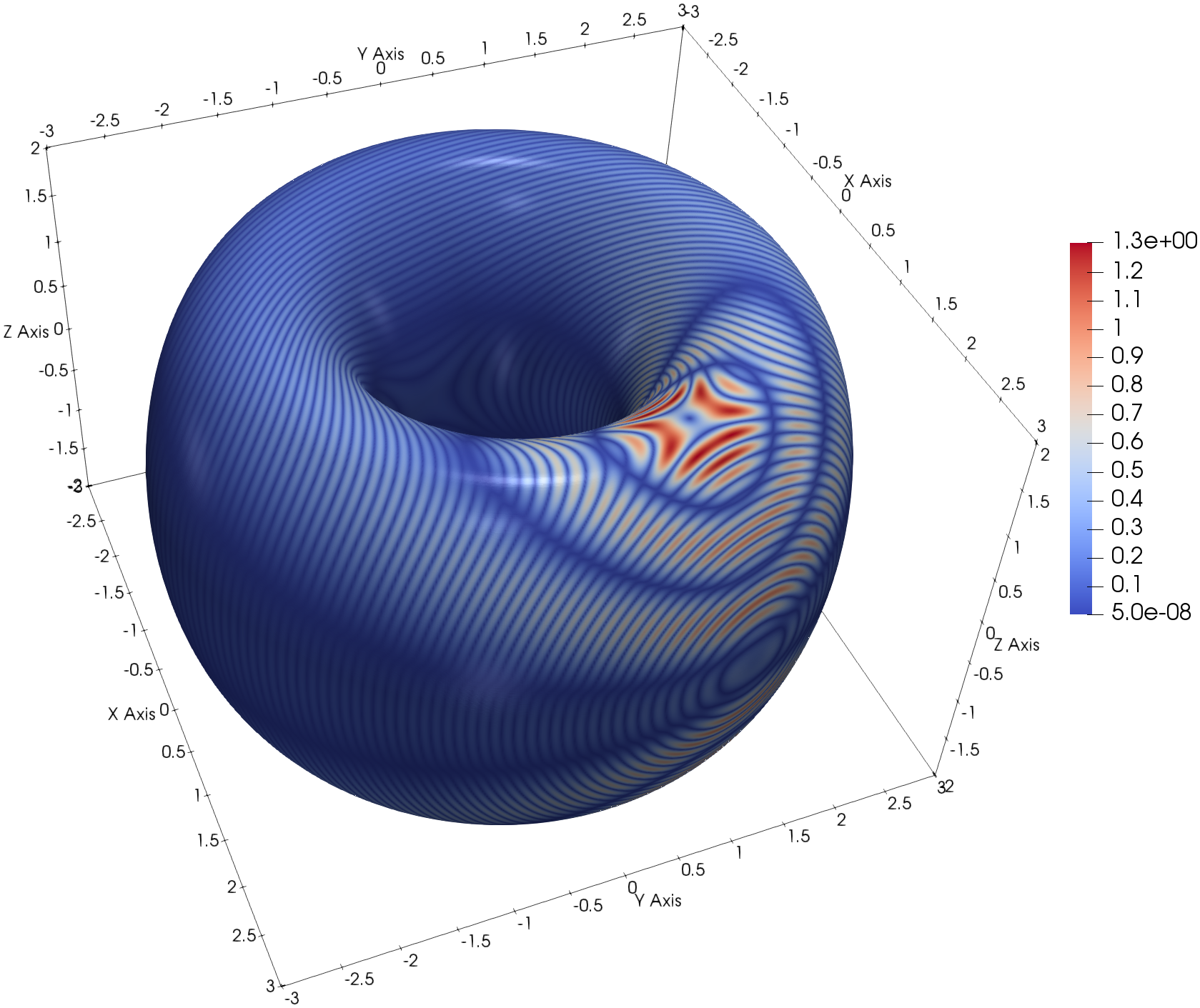}
    \caption{The test geometry and data for $\lambda=0.12$.}
  \end{subfigure}
  \caption{The geometry for our convergence test: a 1-2 torus
    centered at the origin. Shown is $|\realpart \bEin|$ along the surface
    of the scatterer.}
  \label{fig_testgeo}
\end{figure}

\begin{table}[!b]
  \begin{center}
    \caption{Convergence results of the solver for scattering from a
      perfect electric conductor.  In each case,
      $\vct{c}_0 = (0.19, 0.67, 1.00)$, $r_0 = 0.20$, and the integral
      equation was discretized using 8th-order hybrid
      Gauss-Trapezoidal quadrature rules~\cite{alpert}.}
    \label{tab_conv_pec}
    
    \begin{subtable}{.9\linewidth}
      \centering
      \caption{Geometry given in Figure~\ref{fig_testgeo} with $\lambda =
        6.0$.}
      \label{tab_conv_pec1}
      \begin{tabular}{|c|c|c|c|c|c|c|c|} \hline
        $N$  &  $L$ & $NL$ &  $\Res\Gamma$ & $t$ & $t_{G_m}$ & Err($\bE$)
        & Err($\bH$) \\ \hline\hline
      65  &  64 & 4,160 & 3.1E-05 & 0.6 & 0.04 &1.1E-07 & 1.1E-07 \\ \hline
      129 & 128 & 16,512 & 2.2E-09 & 1.6 & 0.17  &1.2E-09 & 1.1E-09 \\ \hline
      257 & 256 & 65,792 & 6.4E-15 & 12.8 & 0.57  &1.6E-12 & 1.6E-12 \\ \hline
      513 & 512 & 262,656 & 3.9E-15 & 157.3 & 1.85  &1.0E-12 & 1.0E-12 \\ \hline
      \end{tabular}
    \end{subtable}\\
    \vspace{\baselineskip}
    \begin{subtable}{.9\linewidth}
      \centering
      \caption{Geometry given in Figure~\ref{fig_testgeo} with $\lambda =
        0.12$.}
      \label{tab_conv_pec2}
      \begin{tabular}{|c|c|c|c|c|c|c|c|} \hline
        $N$  &  $L$ & $NL$ & $\Res\Gamma $ & $t$ & $t_{G_m}$ & Err($\bE$)
        & Err($\bH$) \\ \hline\hline
        65  &  64 & 4,160 & 3.1E-05 & 0.6  & 0.08 & 1.0E+00  & 1.0E+00  \\ \hline
        129 & 128 & 16,512 & 2.2E-09 & 1.8 & 0.25 & 1.0E+00 & 1.0E+00 \\ \hline
        257 & 256 & 65,792 & 6.4E-15 & 18.3 & 0.80 & 3.1E-01 & 1.8E-01 \\ \hline
        513 & 512 & 262,656 & 3.9E-15 & 415.5 & 3.09  & 6.1E-06 & 6.1E-06 \\ \hline
      \end{tabular}
    \end{subtable}
  \end{center}
\end{table}
\clearpage
}

The geometry used in the convergence studies for the PEC case is given
by:
\begin{equation}
  \begin{aligned}
    x(s,t) &= \left( 2 + \cos t \right) \cos s \\
    y(s,t) &= \left( 2 + \cos t \right) \sin s\\
    z(s,t) &= 2\, \sin t 
  \end{aligned}
\end{equation}
for $s,t \in [0,2\pi)$.
The geometry used in the convergence studies for the dielectric case is given
by:
\begin{equation}
  \begin{aligned}
    x(s,t) &= \left( 3 + 2.5 \, \cos t \right) \cos s\\
    y(s,t) &= \left( 3 + 2.5 \, \cos t \right) \sin s\\
    z(s,t) &= \sin t
  \end{aligned}
\end{equation}
for $s,t \in [0,2\pi)$.

In the perfectly conducting examples, the error was computed on a
sphere of radius 5 centered at the origin (and therefore exterior to, and
enclosing the scatterer).
As we can see, for the low-frequency perfectly conducting 
example in Table~\ref{tab_conv_pec1}, the solver rapidly 
converges to
near machine precision (up to the conditioning of the linear
system). For higher frequency problems, as given in
Table~\ref{tab_conv_pec2}
 the solver does not obtain any
accuracy until the discretization is sufficiently
fine to resolve the data, at which point the solver begins to converge
rapidly. The time taken to construct the discretized
integral equation also increases for smaller wavelengths as  more Fourier
modes in the azimuthal direction must be used to resolve the
data. Lastly, note that the time required to evaluate the modal
Green's functions is almost negligible, and scales linearly with the
number of matrix entries. Most of the time spent in the
solver is in BLAS3 matrix-matrix multiplication.
Similar results are obtained in the case of the dielectric solver in
Table~\ref{tab_conv_die}. In this case, 
the error in the exterior was computed on a
sphere of radius 8 centered at the origin, and the error in the
interior was computed on a sphere of radius 0.15 centered at $(0,3,0)$.

\begin{table}[!b]
   \begin{center}
      \caption{Convergence results of the solver for scattering from a
        dielectric region.  The exterior test field was generated from
        a current loop located at $(2.83, 1.00, 0.10)$ with radius
        $0.20$ and the interior test field was generated from a
        current loop located at $(-2.83, -1.00, 2.55)$ with radius
        $0.20$.  The integral equation was discretized using 8th-order
        hybrid Gauss-Trapezoidal quadrature rules~\cite{alpert}. The
        geometry and data are shown in Figure~\ref{fig_testgeo2}.}
      \label{tab_conv_die}
      \begin{tabular}{|c|c|c|c|c|c|c|c|c|} \hline
        $N$  &  $L$ & $NL$ & $\Res\Gamma$ & $t$ & Err($\bE_0$)
        & Err($\bH_0$)
        & Err($\bE_1$)
        & Err($\bH_1$) \\ \hline \hline
        65  &  64 & 4,160 & 6.2E-04 & 0.9  & 7.0E-02 & 2.2E-02 
        & 1.4E-04 & 5.3E-05\\ \hline
        129 & 128 & 16,512 & 1.1E-06 & 8.3  & 5.4E-06 & 5.2E-06 
        & 3.0E-07 & 1.3E-07 \\ \hline
        257 & 256 & 65,792 & 5.7E-12 & 85.7 & 3.6E-09 & 3.6E-09 
        & 2.4E-10 & 1.5E-09 \\ \hline
        513 & 512 & 262,656 & 4.1E-15 & 768.5 & 1.3E-11 & 1.4E-11 
        & 9.4E-12 & 3.5E-12 \\ \hline
      \end{tabular}
    \end{center}
\end{table}

\subsection{Varying the order of accuracy of the quadrature rule}

\begin{figure}[t]
  \centering
  \begin{subfigure}[b]{.45\linewidth}
    \centering
    \includegraphics[width=.95\linewidth]{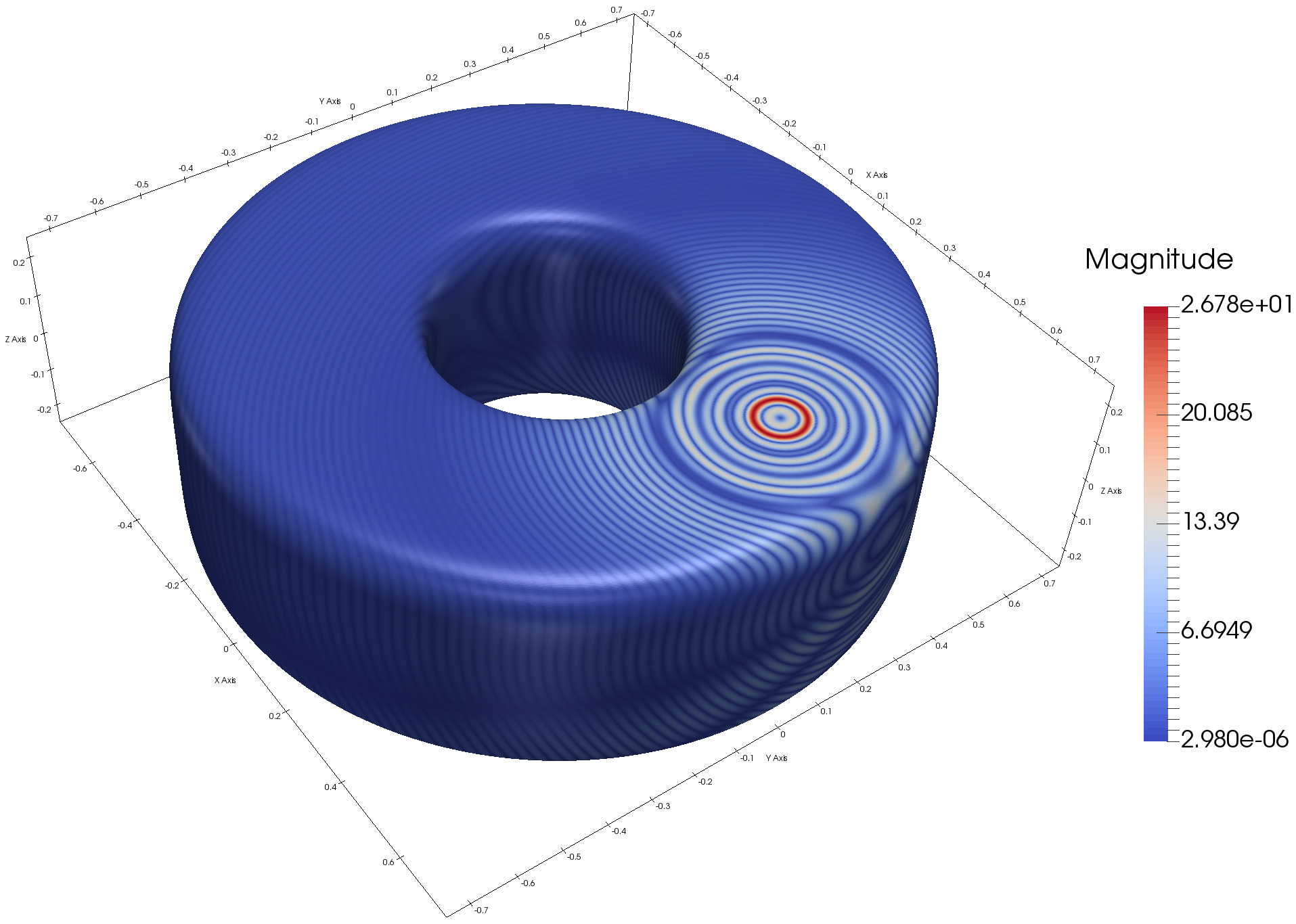}
    \caption{Super-ellipse test geometry (a thick washer)
      and data~$|\realpart \bEin|$ with~$b=0.25$ and~$\lambda=1/32$.}
    \label{fig_pipe1}
  \end{subfigure}
  \hfill
  \begin{subfigure}[b]{.45\linewidth}
    \centering
    \includegraphics[width=.95\linewidth]{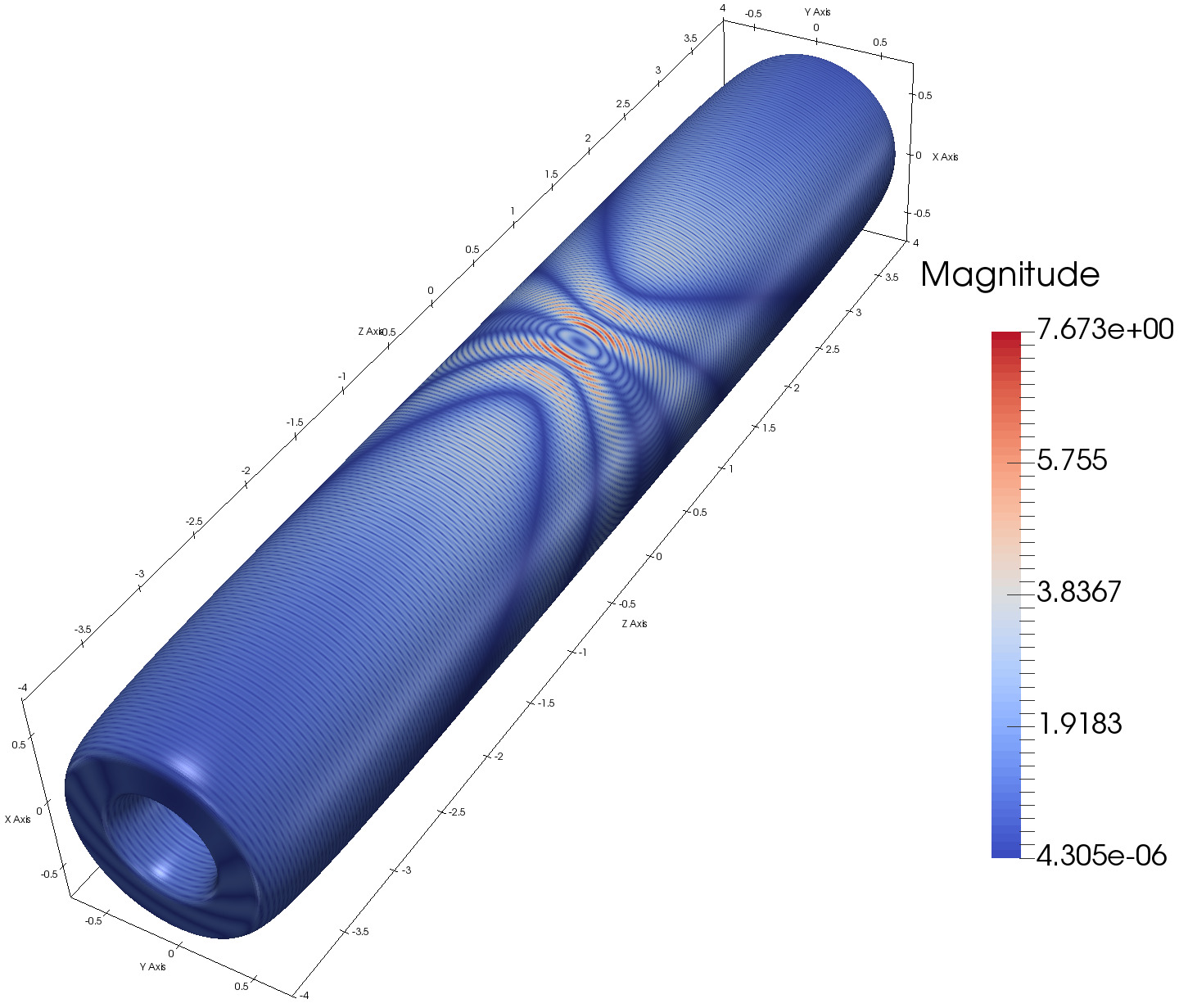}
    \caption{Super-ellipse test geometry (a truncated pipe) and
      data~$|\realpart \bEin|$ with $b=4.0$ and $\lambda=0.06$.}
    \label{fig_pipe2}
  \end{subfigure}
  \caption{The test geometries for demonstrating convergence as a
    function of order of the discretization.}
  \label{fig_pipe}
\end{figure}

We now investigate the effect of the order of accuracy of the quadrature rule on 
the achieved precision, while also considering some
more complicated geometries. We will use the scatterers depicted
in Figure~\ref{fig_pipe}, whose surfaces are given by
\begin{equation}
\left( \frac{r - r_0}{a} \right)^p + \left( \frac{z - z_0}{b}
\right)^p = 1 .
\end{equation}
The generating curve (i.e. a slice of the above surface for fixed
$\theta$) is a \emph{super-ellipse} in
the $\theta z$-plane. The parameter $p$ controls the curvature at the
\emph{corners}, and $a$ and $b$ control the aspect ratio. Geometries
such as this are notoriously hard to resolve using equispaced sampling
schemes. This
generating curve 
can be parameterized with respect to central angle
relative to~$(r_0,0,z_0)$ as:
\begin{equation}
\vct{\gamma}(t) = \left( 
r_0 + \frac{\cos t}{R(t)} \right) \rhat + 
\left( 
  z_0 + \frac{\sin t}{R(t)} \right) \khat,
\end{equation}
where
\begin{equation}
R(t) =  \left( \left(\frac{\cos t}{a} \right)^p
+ \left(\frac{\sin t}{b} \right)^p \right)^{1/p}
\end{equation}
is the distance from the point $(r_0,0,z_0)$ to the boundary (in some
$\theta$-slice).  This is not an arclength parameterization.
In the following calculations, we
have set~$p=6$, $(r_0,0,z_0) = (0.5,0,0)$, $a=0.25$, and the
curve has been discretized (i.e. resampled) in arclength. The height
parameter~$b$ is varied between experiments and given in the captions.

The column labelled $\lambda$ denotes the wavelength,
$\lambda = k/2\pi$.  As before, the columns labelled~$N$ and~$L$
denote the number discretization points along the generating curve, in
the azimuthal direction, respectively. The columns
labelled~$N/\lambda$,~$L/\lambda$, and $D/\lambda$ denote the number
of discretization points {\em per wavelength} along the generating curve,
the number of discretization points {\em per wavelength} along the maximum
circumference of the object in the azimuthal direction, and the
diameter of the smallest sphere enclosing the object in wavelengths,
respectively.  The column labelled~$q$ denotes the order of accuracy of the
quadrature used~\cite{alpert} and $t$ indicates the solution time in seconds.

We see in Table~\ref{tab_pipe1} similar accuracies in the solver at
large wavelengths, independent of the quadrature order.
However, as the wavelength is decreased and the
discretization is refined (in order to maintain a constant number
of points per wavelength), the 2nd order scheme ceases to converge past a few
digits. This is a well-known phenomenon encountered in high
frequency wave propagation problems: in order to maintain precision,
the order of accuracy of the discretization must be increased as the wavelength
decreases (assuming, of course, that the number of discretization points per
wavelength is held constant).

Results for a similar test on an elongated super-ellipse are contained
in Table~\ref{tab_pipe2}.
At higher frequencies, this pipe-like geometry
is particularly difficult to handle, as the interior is many cubic wavelengths
in size and is a nearly resonant cavity.
This example illustrates the importance of
discretizing such problems using high-order
methods.
Table~\ref{tab_pipe2} shows accuracy results for various
order quadrature rules 
at a fixed discretization size for a relatively
small wavelength. With no extra computation time, many more digits of
accuracy can be obtained using a high-order discretization than a
low-order one.
In both examples the error in the computed solution was
determined at 50 points on a sphere of radius 8 centered at the origin.


\begin{table}[!b]
  \begin{center}
    \caption{Accuracy results of the solver for PEC scattering from the
      super-ellipse geometry at relatively high frequencies.  The
      exterior test field was generated from a current loop located at
      $(0.47, 0.17, 0.10)$ with radius $0.1$. The integral equation was
      discretized using varying orders of hybrid Gauss-Trapezoidal
      quadrature rules~\cite{alpert}. The geometries and data are shown
      in Figure~\ref{fig_pipe}.}
    \label{tab_pipe}
    \begin{subtable}{.95\linewidth}
      \centering
      \caption{Convergence for super-ellipse 1 in
        Figure~\ref{fig_pipe1}, maintaining a constant number of points
        per wavelength 
        along the generating curve ($N/\lambda$) and along the maximum circumference
        of the object in the
        azimuthal direction ($L/\lambda$).}
      \label{tab_pipe1}
      \resizebox{\linewidth}{!}{
      \begin{tabular}{|c|c|c|c|c|c|c|c|c|c|c|c|} \hline
        $\lambda$ & $N$  & $N/\lambda$ & $L$ & $L/\lambda$ & $NL$
            & $D/\lambda$ & $\Res\Gamma$
            & $q$ & $t$ & Err($\bE$) & Err($\bH$) \\ \hline\hline
        1/2 & 41 & 11.2 & 100 & 10.6  & 4,100 & 3.0
            & 1.2E-02 & 2 & 0.4 & 4.7E-03 & 4.8E-03 \\ \hline
        1/4 & 81 & 11.1 & 190 & 10.0 & 15,390 & 6.1
            & 3.5E-04 & 2 & 0.6 & 1.4E-03 & 1.4E-03 \\ \hline
        1/8 & 151 & 10.3  & 380 &  10.0 & 57,380
            &  12.2 & 1.2E-06 & 2 & 2.1 & 5.9E-04 & 5.9E-04 \\ \hline
        1/16 & 301 & 10.3 & 760 &  10.0 & 228,760
            &  24.5 & 8.4E-12 & 2 & 10.0 & 4.3E-04 & 4.3E-04 \\ \hline
        1/32 & 591 & 10.0 & 1510 &  10.0 & 892,410
            &  49.0 & 6.8E-15 & 2 & 68.4 & 4.5E-04 & 4.5E-04 \\ \hline\hline
        1/2 & 41 & 11.2 & 100 & 10.6  & 4,100 & 3.0
            & 1.2E-02 & 4 & 0.2 & 2.4E-03 & 2.6E-03 \\ \hline
        1/4 & 81 & 11.1 & 190 & 10.0 & 15,390 & 6.1
            & 3.5E-04 & 4 & 0.6 & 2.3E-04 & 2.3E-04 \\ \hline
        1/8 & 151 & 10.3  & 380 &  10.0 & 57,380
            &  12.2 & 1.2E-06 & 4 & 2.0 & 3.8E-05 & 3.8E-05 \\ \hline
        1/16 & 301 & 10.3 & 760 &  10.0 & 228,760
            &  24.5 & 8.4E-12 & 4 & 10.0 &1.0E-05 & 1.0E-05 \\ \hline
        1/32 & 591 & 10.0 & 1510 &  10.0 & 892,410
            &  49.0 & 6.8E-15 & 4 & 67.7 & 1.0E-07 & 1.0E-07 \\ \hline \hline
        1/2 & 41 & 11.2 & 100 & 0.6  & 4,100 & 3.0
            & 1.2E-02 & 8 & 0.2 & 1.9E-03 & 2.1E-03 \\ \hline
        1/4 & 81 & 11.1 & 190 & 10.0 & 15,390 & 6.1
            & 3.5E-04 & 8 & 0.6 & 1.9E-05 & 2.0E-05 \\ \hline
        1/8 & 151 & 10.3  & 380 &  10.0 & 57,380
            &  12.2 & 1.2E-06 & 8 & 1.9 & 9.1E-07 & 9.1E-07 \\ \hline
        1/16 & 301 & 10.3 & 760 &  10.0 & 228,760
            &  24.5 & 8.4E-12 & 8 & 10.2 & 4.6E-08 & 4.6E-08 \\ \hline
        1/32 & 591 & 10.0 & 1510 &  10.0 & 892,410
            &  49.0 & 6.8E-15 & 8 & 67.8 & 3.0E-08 & 3.0E-08 \\ \hline \hline
        1/2 & 41 & 11.2 & 100 & 10.6  & 4,100 & 3.0
            & 1.2E-02 & 16 & 0.2 & 1.9E-03 & 2.0E-03 \\ \hline
        1/4 & 81 & 11.1 & 190 & 10.0 & 15,390 & 6.1
            & 3.5E-04 & 16 & 0.7 & 6.2E-06 & 7.1E-06 \\ \hline
        1/8 & 151 & 10.3  & 380 &  10.0 & 57,380
            &  12.2 & 1.2E-06 & 16 & 2.2 & 1.6E-08 & 1.6E-08 \\ \hline
        1/16 & 301 & 10.3 & 760 &  10.0 & 228,760
            &  24.5 & 8.4E-12 & 16 & 10.5 & 4.5E-10 & 4.5E-10 \\ \hline
        1/32 & 591 & 10.0 & 1510 &  10.0 & 892,410
            &  49.0 & 6.8E-15 & 16 & 68.5 & 3.8E-10 & 3.8E-10 \\ \hline
    \end{tabular}}
  \end{subtable}\\
  \vspace{\baselineskip}
\begin{subtable}{.9\linewidth}
  \centering
    \caption{Accuracy obtained for varying quadrature orders $q$ in the
      discretization of
      the elongated super-ellipse (truncated pipe) shown
      given in Figure~\ref{fig_pipe2} with $\lambda =
    0.06$. The pipe is approximately 135$\lambda$ long,
    making the scattering problem rather ill-conditioned.}
  \label{tab_pipe2}
      \resizebox{\linewidth}{!}{
      \begin{tabular}{|c|c|c|c|c|c|c|c|c|c|c|c|} \hline
        $\lambda$ & $N$  & $N/\lambda$ & $L$ & $L/\lambda$ & $NL$
            & $D/\lambda$ & $\Res\Gamma$
            & $q$ & $t$ & Err($\bE$) & Err($\bH$) \\
        \hline\hline
        0.06 & 1381 & 5.0 & 395 & 5.0 & 545,495 & 134.7
            & 3.2E-07 & 2 &  483.5 & 1.3E-02 & 1.3E-02 \\ \hline
        0.06 & 1381 & 5.0 & 395 & 5.0 & 545,495 & 134.7
            & 3.2E-07 & 4 & 504.2 & 1.2E-03 & 1.2E-03 \\ \hline
        0.06 & 1381 & 5.0 & 395 & 5.0 & 545,495 & 134.7
            & 3.2E-07 & 8 & 491.5 & 1.9E-04 & 1.9E-04 \\ \hline
        0.06 & 1381 & 5.0 & 395 & 5.0 & 545,495 & 134.7
            & 3.2E-07 & 16 & 486.6 & 1.1E-06 & 1.1E-06 \\ \hline
      \end{tabular}
      }
    \end{subtable}
  \end{center}
\end{table}


\subsection{Computation of radar cross-section}

With the previous high-order convergence results in mind, the canonical
application of integral equation-based solvers for electromagnetic
scattering phenomena is the computation of the \emph{radar
  cross-section} (RCS) of an object. Essentially, the radar cross
section characterizes the far-field scattering response of an object
to an incoming electromagnetic plane-wave and is critical in 
inverse scattering and object recognition applications.

\begin{figure}[!t]
  \centering
  \begin{subfigure}[b]{.4\linewidth}
    \centering
    \includegraphics[width=.95\linewidth]{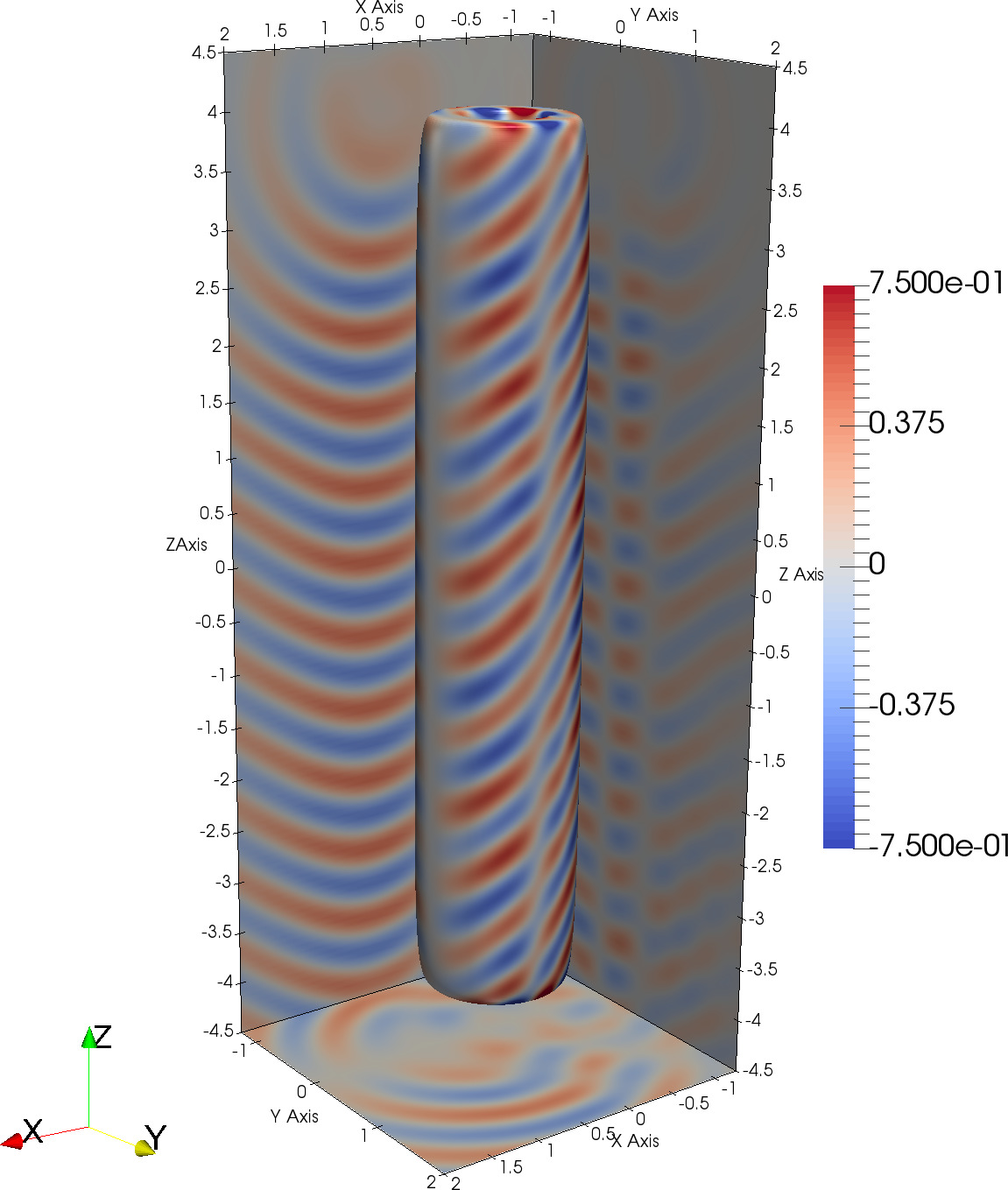}
    \caption{The induced generalized Debye source~$\rho$
      and the total field 
      on an enclosing box.}
    \label{fig_mrcs1}
  \end{subfigure}
  \hfill
  \begin{subfigure}[b]{.55\linewidth}
    \centering
    \includegraphics[width=.95\linewidth]{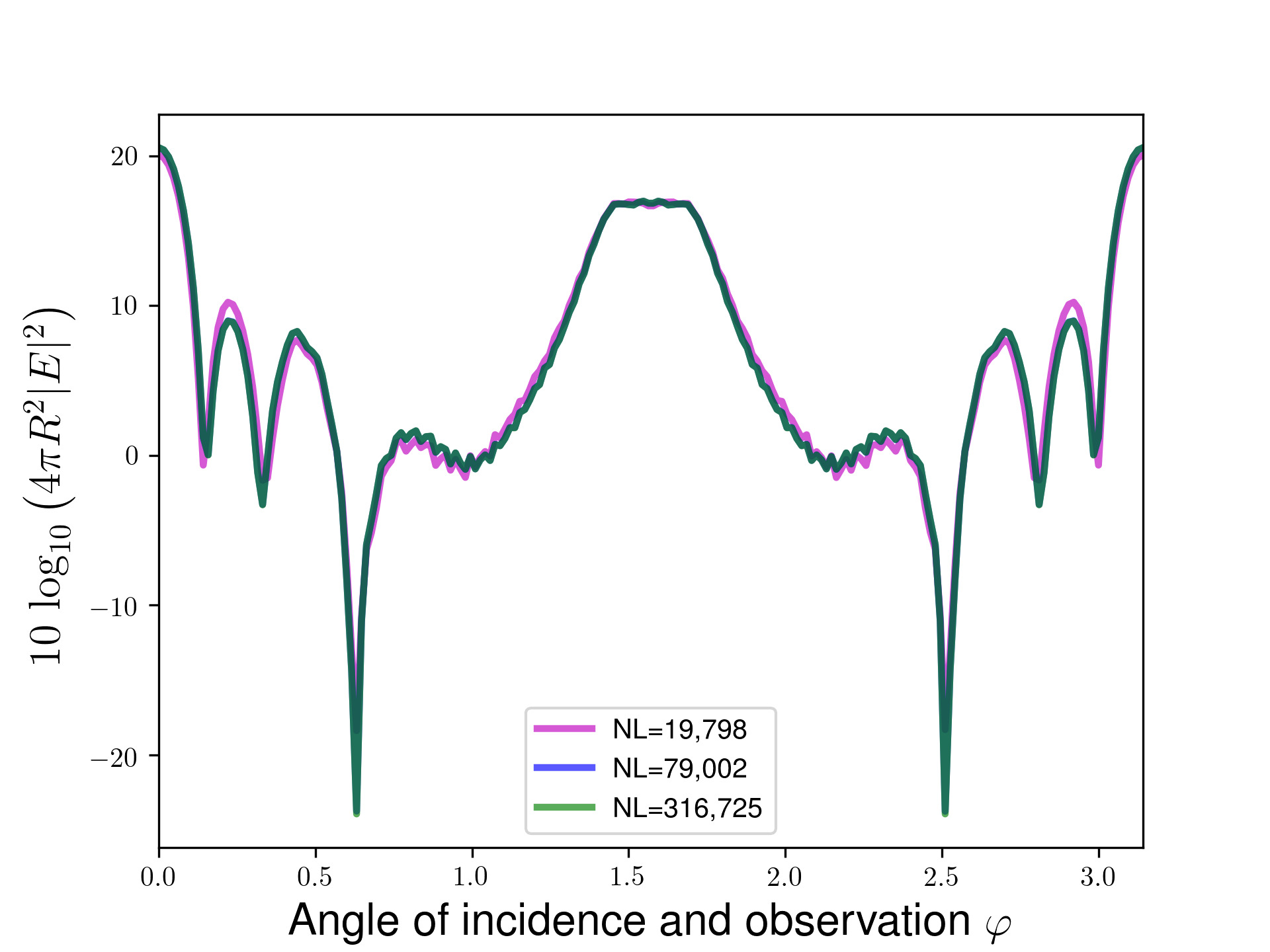}
    \caption{The monostatic radar cross section of the pipe geometry
      for a horizontally polarized incoming plane wave.}
    \label{fig_mrcs2}
  \end{subfigure}
  \caption{Plots of the scattered field from an incident plane wave and the
    monostatic radar cross section for the pipe geometry.}
  \label{fig_mrcs}
\end{figure}

An incoming electromagnetic plane-wave propagating in the direction
$\vct{u}$ with polarization~$\vct{p}$ is given by
\begin{equation}
  \begin{aligned}
    \bEin(\bx) &= \vct{p}  \, e^{ik  {\vec{u}} \cdot \bx}, \\
    \bHin(\bx) &= \vec{u} \times \vct{p} \,
        e^{ik  \vec{u} \cdot \bx},
  \end{aligned}
\end{equation}
where it is assumed that $\vct{u}$ and $\vct{p}$ are unit vectors. In
order for Maxwell's equations to be satisfied, we must have
$\vct{u} \cdot \vct{p} = 0$. The \emph{monostatic RCS} (MRCS) of an
object~$\Omega$ degscribes the backscatter of an object in the
direction of propagation opposite that of the incoming plane wave. 
For a body of revolution, it is standard to compute
only the $\phi$-dependence (polar angle) of the MRCS,
since for fixed $\pm z$-polarizations
the~$\theta$-dependence is constant (due to rotational symmetry).
Letting~$(\rho,\phi,\theta)$
denote the spherical coordinates with~$\theta$ the azimuthal
angle, if we set the direction of propagation of the incoming
plane-wave to be $\vct{u} = -\sin\phi \, \ihat -\cos\phi \, \khat$,
then a horizontally-polarized plane wave is given by:
\begin{equation}
  \begin{aligned}
    \bEin(\bx) &= \jhat  \, e^{-ik(x\sin\phi + z\cos\phi)}, \\
    \bHin(\bx) &= (\cos\phi \, \ihat - \sin\phi \, \khat)  \,
        e^{-ik(x\sin\phi + z\cos\phi)},
  \end{aligned}
\end{equation}
where as before, $(\ihat, \jhat, \khat)$ are the unit vectors in the
Cartesian coordinate system.
As a function of $\phi$, the MRCS in this setup is then
approximately~\cite{jackson,jin-2010}:
\begin{equation}
  \begin{aligned}
    \text{MRCS} (\phi) &\approx 4\pi R^2 \,
    \frac{|\bE(\bx_R(\phi)) |^2}{| \bEin(\bx_R(\phi)) |^2} \\
    &= 4\pi R^2 \, |\bE(\bx_R(\phi)) |^2
  \end{aligned}
\end{equation}
since~$| \bEin | = 1$, where $R$ is some fixed distance (set to be 10
in the subsequent experiments) from the
center of the object (assumed to be the origin of the
coordinate system in this example) and
$\bx_R(\phi) = R\sin\phi \, \ihat + R\cos\phi \, \khat$.  Take care to note
that the observation point~$\bx_R$ lies in the direction
\emph{opposite} the direction of propagation of the incoming plane
wave.  Shown in Figure~\ref{fig_mrcs1} is a plot of the real-part of
the induced generalized Deybe source ~$\rho$ along~$\Omega$,
$\realpart \rho$ (see equation~\eqref{eq_jbuild}), and the real-part
of the $x$-component of the total field, $\realpart E^{\text{tot}}_x$,
along the enclosing walls. The incoming plane wave
was assumed to be propagating in the direction~$(-1/\sqrt{2},0,-1/\sqrt{2})$ with
polarization~$(0,1,0)$. In Figure~\ref{fig_mrcs2}
we plot the MCRS as a function of the polar angle of incidence, $\phi
\in [0,\pi]$, for various refinements of the discretization. 
Following convention, we
report the MRCS on a decibel scale: $10 \, \log_{10} \text{MRCS}$.
In these examples, we set $\lambda= 0.5$ so that the pipe is
approximately 16 wavelengths long. The pipe was discretized with
$N=129$, $N=513$, and $N=1025$ points along the generating curve, 
and $L$ was set so as to maintain the same
number of points per wavelength along the azimuthal direction.

These experiments correspond to true scattering problems, for which the
exact solution is not known. Since the data are very regular and smooth
(in fact smoother than the data in the previous examples) the accuracy
can be estimated from Table~\ref{tab_pipe2}.  The accuracy can also be
estimated by a self-convergence study, as shown in the column
labeled~\emph{Err(MRCS)} in Table~\ref{tab_rcs}. This
error is the relative $\ell^2$ error in the quantity
$10 \, \log_{10} MRCS$.  We see that the errors are approximately
commensurate with the geometry resolution, and comparable to the
results in the high-frequency case in the previous subsection.
Furthermore, since we are using a direct solver, the system can be set
up and factored once, so that each subsequent MRCS calculation merely
requires the application of the inverse to a new right hand side,
greatly accelerating the calculation.

\begin{table}[!b]
  \begin{center}
    \caption{Data for the monostatic radar cross section plot in
      Figure~\ref{fig_mrcs2}. The relative $\ell^2$ error is computed via a
      self-convergence test at 200 equispaced
      points in the polar angle direction at a distance of $R=10$.}
    \label{tab_rcs}
      \begin{tabular}{|c|c|c|c|c|c|c|c|c|} \hline
        $\lambda$ & $N$  & $N/\lambda$ & $L$ & $L/\lambda$ & $NL$
        & $D/\lambda$ & $\Res\Gamma$
        & Err(MRCS) \\ \hline\hline
        0.5 & 257 & 7.8 & 77 & 8.2  & 19,789 & 8.1
                      & 1.2E-02 &  1.1E-01  \\ \hline
        0.5 &  513 & 15.5 & 154 & 16.3  & 79,002 & 8.1 
                      & 1.2E-03  & 2.6E-03  \\ \hline
        0.5 & 1025  & 31.0 & 309 & 32.8  & 316,725 & 8.1
                      & 1.4E-05 &   \\ \hline
      \end{tabular}
  \end{center}
\end{table}

\section{Conclusions}
\label{sec_conclusions}

We have presented a direct solver based on the generalized Debye
integral equation framework for full electromagnetic scattering from
perfect electric conductors or dielectric bodies which are rotationally
symmetric. Unlike most widely used formulations, our approach
is invertible for all passive materials at any frequency, and 
immune from both dense-mesh breakdown and topological low frequency
breakdown.

The scheme makes use of separation of variables
in the azimuthal direction, leading to a collection of smaller,
uncoupled integral equations on curves.
The implementation described in this work relies on a trapezoidal-based
Nystr\"om discretization, restricting its applicability to
geometries which are reasonably well-discretized using equispaced
nodes in some parameterization. Alternative discretizations and quadrature
rules are required for scatterers
whose generating curves contain corner singularities, but the overall approach
would remain the same.
For large scale problems, the bulk of the
computation is spent in BLAS3 routines, assembling the system matrix
(due to the composition of several integral operators in the representation).
These matrix assembly steps are highly parallelizable.  

We are presently working on 
implementing solvers for
the generalized Debye integral equations on more general 
(closed) surfaces in three dimensions.

\section*{Acknowledgments}
We would like to acknowledge Alex Barnett, Johan Helsing, and Gunnar
Martinsson for many useful conversations.  We also gratefully
acknowledge the support of the NVIDIA Corporation with the donation of
a Quadro P6000, used for some of the visualizations presented in this
research.

\bibliographystyle{abbrvnat}
\bibliography{master}

\begin{thebibliography}{52}
\providecommand{\natexlab}[1]{#1}
\providecommand{\url}[1]{\texttt{#1}}
\expandafter\ifx\csname urlstyle\endcsname\relax
  \providecommand{\doi}[1]{doi: #1}\else
  \providecommand{\doi}{doi: \begingroup \urlstyle{rm}\Url}\fi

\bibitem[Ahrens et~al.(2005)Ahrens, Geveci, and Law]{paraview2005}
J.~Ahrens, B.~Geveci, and C.~Law.
\newblock Paraview: {A}n end-user tool for large data visualization.
\newblock In C.~D. Hansen and C.~R. Johnson, editors, \emph{{The Visualization
  Handbook}}. Elsevier, 2005.

\bibitem[Alpert(1999)]{alpert}
B.~Alpert.
\newblock Hybrid {G}auss-trapezoidal quadrature rules.
\newblock \emph{{SIAM} {J}. {S}ci. {C}omput.}, 20\penalty0 (5):\penalty0
  1551--1584, 1999.

\bibitem[Atkinson(1997)]{atkinson_1997}
K.~E. Atkinson.
\newblock \emph{{The Numerical Solution of Integral Equations of the Second
  Kind}}.
\newblock Cambridge University Press, New York, NY, 1997.

\bibitem[Barnett and Betcke(2008)]{barnett2008stability}
A.~H. Barnett and T.~Betcke.
\newblock Stability and convergence of the method of fundamental solutions for
  {H}elmholtz problems on analytic domains.
\newblock \emph{Journal of Computational Physics}, 227\penalty0 (14):\penalty0
  7003--7026, 2008.

\bibitem[Bremer and Gimbutas(2013)]{bremer_2013}
J.~Bremer and Z.~Gimbutas.
\newblock On the numerical evaluation of singular integrals of scattering
  theory.
\newblock \emph{J. Comput. Phys.}, 251:\penalty0 327--343, 2013.

\bibitem[Briggs and Henson(1995)]{briggs_1995}
W.~L. Briggs and V.~E. Henson.
\newblock \emph{{The DFT: An Owner's Manual for the Discrete Fourier
  Transform}}.
\newblock SIAM, Philadelphia, PA, 1995.

\bibitem[Chew et~al.(2001)Chew, Michielssen, Song, and Jin]{Chewfastbook}
W.~C. Chew, E.~Michielssen, J.~M. Song, and J.~M. Jin.
\newblock \emph{Fast and Efficient Algorithms in Computational
  Electromagnetics}.
\newblock Artech House, Inc., Norwood, MA, 2001.

\bibitem[Cohl and Tohline(1999)]{cohl_1999}
H.~S. Cohl and J.~E. Tohline.
\newblock A compact cylindrical {G}reen's function expansion for the solution
  of potential problems.
\newblock \emph{Astrophys. J.}, 527\penalty0 (1):\penalty0 86--101, 1999.

\bibitem[Colton and Kress(1983)]{colton_kress}
D.~Colton and R.~Kress.
\newblock \emph{Integral {E}quation {M}ethods in {S}cattering {T}heory}.
\newblock John Wiley \& Sons, Inc., 1983.

\bibitem[Contopanagos et~al.(2002)Contopanagos, Dembart, Epton, Ottusch,
  Rokhlin, Visher, and Wandzura]{contopanagos-2002}
H.~Contopanagos, B.~Dembart, M.~Epton, J.~J. Ottusch, V.~Rokhlin, J.~L. Visher,
  and S.~M. Wandzura.
\newblock Well-conditioned boundary integral equations for three-dimensional
  electromagnetic scattering.
\newblock \emph{IEEE Trans. Antennas Propag.}, 50\penalty0 (12):\penalty0
  1824--1830, 2002.

\bibitem[Conway and Cohl(2010)]{conway_cohl}
J.~T. Conway and H.~S. Cohl.
\newblock Exact {F}ourier expansion in cylindrical coordinates for the
  three-dimensional {H}elmholtz {G}reen function.
\newblock \emph{Z. {A}ngew. {M}ath. {P}hys.}, 61:\penalty0 425--442, 2010.

\bibitem[Cools et~al.(2009)Cools, Andriulli, Olyslager, and
  Michielssen]{cools-2009}
K.~Cools, F.~P. Andriulli, F.~Olyslager, and E.~Michielssen.
\newblock {Nullspaces of MFIE and Calderon Preconditioned EFIE Operators
  Applied to Toroidal Surfaces}.
\newblock \emph{IEEE Trans. Antennas Propag.}, 57\penalty0 (10):\penalty0
  3205--3215, 2009.

\bibitem[Epstein and Greengard(2010)]{EpGr}
C.~L. Epstein and L.~Greengard.
\newblock Debye sources and the numerical solution of the time harmonic
  {M}axwell equations.
\newblock \emph{Commun. Pure Appl. Math.}, 63\penalty0 (4):\penalty0 413--463,
  2010.

\bibitem[Epstein et~al.(2013{\natexlab{a}})Epstein, Gimbutas, Greengard,
  Kl\"ockner, and O'Neil]{epstein-2013}
C.~L. Epstein, Z.~Gimbutas, L.~Greengard, A.~Kl\"ockner, and M.~O'Neil.
\newblock A consistency condition for the vector potential in
  multiply-connected domains.
\newblock \emph{IEEE Trans. Magn.}, 49\penalty0 (3):\penalty0 1072--1076,
  2013{\natexlab{a}}.

\bibitem[Epstein et~al.(2013{\natexlab{b}})Epstein, Greengard, and
  O'Neil]{EpGrOn}
C.~L. Epstein, L.~Greengard, and M.~O'Neil.
\newblock Debye sources and the numerical solution of the time harmonic
  {M}axwell equations {II}.
\newblock \emph{Commun. Pure. Appl. Math.}, 66\penalty0 (5):\penalty0 753--789,
  2013{\natexlab{b}}.

\bibitem[Epstein et~al.(2015)Epstein, Greengard, and O'Neil]{EpGrOn2}
C.~L. Epstein, L.~Greengard, and M.~O'Neil.
\newblock Debye sources, {B}eltrami fields, and a complex structure on
  {M}axwell fields.
\newblock \emph{Commun. Pure. Appl. Math.}, 68:\penalty0 2237--2280, 2015.

\bibitem[Fairweather and Karageorghis(1998)]{fairweather1998method}
G.~Fairweather and A.~Karageorghis.
\newblock The method of fundamental solutions for elliptic boundary value
  problems.
\newblock \emph{Advances in Computational Mathematics}, 9\penalty0
  (1-2):\penalty0 69, 1998.

\bibitem[Frankel(2011)]{frankel}
T.~Frankel.
\newblock \emph{The Geometry of Physics}.
\newblock Cambridge University Press, New York, NY, 2011.

\bibitem[Frittelli and Sgura(2016)]{frittelli_2016}
M.~Frittelli and I.~Sgura.
\newblock Virtual element method for the {L}aplace-{B}eltrami equation on
  surfaces.
\newblock 2016.
\newblock arXiv:1612.02369 [math.NA].

\bibitem[Gedney and Mittra(1990)]{gedney1990use}
S.~D. Gedney and R.~Mittra.
\newblock The use of the {FFT} for the efficient solution of the problem of
  electromagnetic scattering by a body of revolution.
\newblock \emph{IEEE Trans. Antennas Propag.}, 38\penalty0 (3):\penalty0
  313--322, 1990.

\bibitem[Gil et~al.(2007)Gil, Segura, and Temme]{gil_2007}
A.~Gil, J.~Segura, and N.~M. Temme.
\newblock \emph{{Numerical Methods for Special Functions}}.
\newblock SIAM, Philadelphia, PA, 2007.

\bibitem[Gimbutas and Greengard(2013)]{gimbutas-2013}
Z.~Gimbutas and L.~Greengard.
\newblock {Fast multi-particle scattering: A hybrid solver for the Maxwell
  equations in microstructured materials}.
\newblock \emph{J. Comput. Phys.}, 232\penalty0 (1):\penalty0 22--32, 2013.

\bibitem[Gustafsson(2010)]{gustafsson2010accurate}
M.~Gustafsson.
\newblock Accurate and efficient evaluation of modal {G}reen's functions.
\newblock \emph{J. Electromagnet. Wave.}, 24\penalty0 (10):\penalty0
  1291--1301, 2010.

\bibitem[Hao et~al.(2014)Hao, Barnett, Martinsson, and Young]{hao_2014}
S.~Hao, A.~H. Barnett, P.-G. Martinsson, and P.~Young.
\newblock High-order accurate {N}ystr\"om discretization of integral equations
  with weakly singular kernels on smooth curves in the plane.
\newblock \emph{Adv. Comput. Math.}, 40:\penalty0 245--272, 2014.

\bibitem[Helsing and Karlsson(2014)]{helsing_2014}
J.~Helsing and A.~Karlsson.
\newblock An explicit kernel-split panel-based {N}ystr\"om scheme for integral
  equations on axially symmetric surfaces.
\newblock \emph{J. Comput. Phys.}, 272:\penalty0 686--703, 2014.

\bibitem[Helsing and Karlsson(2015)]{helsing2015determination}
J.~Helsing and A.~Karlsson.
\newblock Determination of normalized magnetic eigenfields in microwave
  cavities.
\newblock \emph{IEEE Transactions on Microwave Theory and Techniques},
  63\penalty0 (5):\penalty0 1457--1467, 2015.

\bibitem[Imbert-Gerard and Greengard(2017)]{imbertgerard_2017}
L.-M. Imbert-Gerard and L.~Greengard.
\newblock Pseudo-spectral methods for the {L}aplace-{B}eltrami equation and the
  {H}odge decomposition on surfaces of genus one.
\newblock \emph{Numer. Methods Partial. Differ. Equ.}, 33\penalty0
  (3):\penalty0 941--955, 2017.

\bibitem[Jackson(1999)]{jackson}
J.~D. Jackson.
\newblock \emph{Classical Electrodynamics}.
\newblock Wiley, New York, NY, 3rd edition, 1999.

\bibitem[Jin(2010)]{jin-2010}
J.-M. Jin.
\newblock \emph{Theory and Computation of Electromagnetic Fields}.
\newblock IEEE Press, Piscataway, NJ, 2010.

\bibitem[Kapur and Long(1998)]{kapur-1998}
S.~Kapur and D.~E. Long.
\newblock {IES3}: Efficient electrostatic and electromagnetic simulation.
\newblock \emph{Comput. Sci. \& Engrg.}, pages 60--67, 1998.

\bibitem[Kirsch and Monk(1995)]{kirsch1995finite}
A.~Kirsch and P.~Monk.
\newblock A finite element/spectral method for approximating the time-harmonic
  {M}axwell system in ${R}^3$.
\newblock \emph{SIAM J. Appl. Math.}, 55\penalty0 (5):\penalty0 1324--1344,
  1995.

\bibitem[Kucharski(2000)]{Kucharski2000}
A.~A. Kucharski.
\newblock {A method of moments solution for electromagnetic scattering by
  inhomogeneous dielectric bodies of revolution}.
\newblock \emph{IEEE Trans. Antennas Propag.}, 48\penalty0 (8):\penalty0
  1202--1210, 2000.

\bibitem[Liu and Barnett(2016)]{liu2016efficient}
Y.~Liu and A.~H. Barnett.
\newblock Efficient numerical solution of acoustic scattering from
  doubly-periodic arrays of axisymmetric objects.
\newblock \emph{J. Comput. Phys.}, 324:\penalty0 226--245, 2016.

\bibitem[Mautz and Harrington(1979)]{mautz1979electromagnetic}
J.~R. Mautz and R.~F. Harrington.
\newblock Electromagnetic scattering from a homogeneous material body of
  revolution.
\newblock \emph{Archiv Elektronik und Uebertragungstechnik}, 33:\penalty0
  71--80, 1979.

\bibitem[Monk(1992)]{monk1992finite}
P.~Monk.
\newblock A finite element method for approximating the time-harmonic {M}axwell
  equations.
\newblock \emph{Numerische Mathematik}, 63\penalty0 (1):\penalty0 243--261,
  1992.

\bibitem[M{\"u}ller(1969)]{muller}
C.~M{\"u}ller.
\newblock \emph{Foundations of the {M}athematical {T}heory of {E}lectromagnetic
  {W}aves}.
\newblock Springer-Verlag, Berlin, Heidelberg, 1969.

\bibitem[Nedelec(2001)]{nedelec}
J.-C. Nedelec.
\newblock \emph{{Acoustic and Electromagnetic Equations}}.
\newblock Springer, New York, NY, 2001.

\bibitem[Olver et~al.(2010)Olver, Lozier, Boisvert, and Clark]{nist}
F.~W. Olver, D.~W. Lozier, R.~F. Boisvert, and C.~W. Clark.
\newblock \emph{NIST Handbook of Mathematical Functions}.
\newblock Cambridge University Press, New York, NY, USA, 1st edition, 2010.
\newblock ISBN 0521140633, 9780521140638.

\bibitem[O'Neil(2018)]{oneil2017}
M.~O'Neil.
\newblock {Second-kind integral equations for the Laplace-Beltrami problem on
  surfaces in three dimensions}.
\newblock \emph{Adv. Comput. Math.}, 2018.
\newblock To appear.

\bibitem[O'Neil and Cerfon(2018)]{oneil2018taylor}
M.~O'Neil and A.~J. Cerfon.
\newblock An integral equation-based numerical solver for {T}aylor states in
  toroidal geometries.
\newblock \emph{J. Comput. Phys.}, 359:\penalty0 263--282, 2018.

\bibitem[Papas(1988)]{papas}
C.~H. Papas.
\newblock \emph{Theory of Electromagnetic Wave Propagation}.
\newblock Dover, New York, NY, 1988.

\bibitem[Sidi and Israeli(1988)]{sidi1988}
A.~Sidi and M.~Israeli.
\newblock Quadrature methods for periodic singular and weakly singular
  {F}redholm integral equations.
\newblock \emph{J. Sci. Comput.}, 3\penalty0 (2):\penalty0 201--231, 1988.

\bibitem[Sifuentes et~al.(2015)Sifuentes, Gimbutas, and
  Greengard]{sifuentes_2015}
J.~Sifuentes, Z.~Gimbutas, and L.~Greengard.
\newblock Randomized methods for rank-deficient linear systems.
\newblock \emph{Elec. Trans. Num. Anal.}, 44:\penalty0 177--188, 2015.

\bibitem[Taskinen and Yla-Oijala(2006)]{taskinen-2006}
M.~Taskinen and P.~Yla-Oijala.
\newblock Current and charge integral equation formulation.
\newblock \emph{IEEE Trans. Antennas Propag.}, 54\penalty0 (1):\penalty0
  58--67, 2006.

\bibitem[Trefethen(2000)]{trefethen_spectral_book}
L.~N. Trefethen.
\newblock \emph{Spectral methods in {MATLAB}}.
\newblock SIAM, Philadelphia, PA, 2000.

\bibitem[Vico et~al.(2013)Vico, Gimbutas, Greengard, and
  Ferrando-Bataller]{vico-2013}
F.~Vico, Z.~Gimbutas, L.~Greengard, and M.~Ferrando-Bataller.
\newblock Overcoming low-frequency breakdown of the magnetic field integral
  equation.
\newblock \emph{IEEE Trans. Antennas Propag.}, 61\penalty0 (3):\penalty0
  1285--1290, 2013.

\bibitem[Vico et~al.(2016)Vico, Ferrando, Greengard, and Gimbutas]{vico_2016}
F.~Vico, M.~Ferrando, L.~Greengard, and Z.~Gimbutas.
\newblock The decoupled potential integral equation for time-harmonic
  electromagnetic scattering.
\newblock \emph{Comm. Pure Appl. Math.}, 69:\penalty0 771--812, 2016.

\bibitem[Viola(1995)]{viola}
M.~S. Viola.
\newblock {A new electric field integral equation for heterogeneous dielectric
  bodies of revolution}.
\newblock \emph{IEEE Transactions on Microwave Theory and Techniques},
  43:\penalty0 230--233, 1995.

\bibitem[Wright et~al.(2015)Wright, Javed, Montanelli, and
  Trefethen]{wright2015}
G.~B. Wright, M.~Javed, H.~Montanelli, and L.~N. Trefethen.
\newblock Extension of {C}hebfun to periodic functions.
\newblock \emph{SIAM J. Sci. Comput.}, 37:\penalty0 C554--C573, 2015.

\bibitem[Yee(1966)]{yee_1966}
K.~Yee.
\newblock Numerical solution of initial boundary value problems involving
  {M}axwell's equations in isotropic media.
\newblock \emph{IEEE Trans. Antennas Propag.}, 14:\penalty0 302--307, 1966.

\bibitem[Young et~al.(2012)Young, Hao, and Martinsson]{young}
P.~Young, S.~Hao, and P.-G. Martinsson.
\newblock A high-order {N}ystr\"om discretization scheme for boundary integral
  equations defined on rotationally symmetric surfaces.
\newblock \emph{J. Comput. Phys.}, 231\penalty0 (11):\penalty0 4142--4159,
  2012.

\bibitem[Youssef(1989)]{youssef1989}
N.~N. Youssef.
\newblock Radar cross section of complex targets.
\newblock \emph{Proc. IEEE}, 77:\penalty0 722--734, 1989.

\end{thebibliography}
\end{document}